\documentclass[10pt]{article}

\usepackage{latexsym, amssymb, amsmath, amsthm, amscd}
\usepackage[all,cmtip]{xy}
\usepackage{amsmath}
\usepackage{amsthm}
\usepackage{amssymb}
\usepackage{relsize}
\usepackage{mathtools}
\usepackage{amsfonts}
\usepackage{amsrefs}
\usepackage{comment}
\usepackage{color,authblk}
\usepackage{tikz}
\usepackage{amscd} 
\usepackage{comment}
\usepackage{mathrsfs}

\usepackage[all]{xy}
\usepackage{stackrel}
\allowdisplaybreaks
\usepackage{graphicx}

\usepackage{authblk}
\usepackage[margin=1in]{geometry}
\usepackage[hidelinks]{hyperref}
\usepackage{color,soul}
\usepackage{fullpage}
\usepackage[all,cmtip]{xy}
\usepackage{amsmath,amscd}
\usepackage{tikz-cd}
\usepackage{tikz}
\usetikzlibrary{positioning}
\usetikzlibrary{shapes.geometric, arrows}
\DeclareMathAlphabet{\matholdcal}{OMS}{cmsy}{m}{n}

\usetikzlibrary{matrix}

\usepackage[mathcal]{euscript}
\usepackage{txfonts}

\usepackage{hyperref}

\numberwithin{equation}{section} \DeclareMathSizes{2}{10}{12}{13}
\newtheorem{thm}{Proposition}[section]
\newtheorem{Thm}[thm]{Theorem}

\newtheorem{cor}[thm]{Corollary}
\newtheorem{lem}[thm]{Lemma}
\newtheorem{defn}[thm]{Definition}
\newtheorem{prop}[thm]{Proposition}
\newtheorem{eg}[thm]{Example}
\numberwithin{thm}{section} 

\title{\Large Gerstenhaber  type  structures on Davydov-Yetter cohomology with coefficients }

\author{Mamta Balodi \footnote{Department of Mathematics, Indian Institute of Technology, Gandhinagar, India. Email: mamta.balodi@gmail.com.}$ \textrm{ }\footnote{MB was supported by WISE-PDF grant DST/WISE-PDF/PM-29/2024 (G) from Department of Science and Technology, Government of India.} \qquad\qquad$ Abhishek Banerjee\footnote{Department of Mathematics, Indian Institute of Science, Bangalore, India. Email: abhishekbanerjee1313@gmail.com.}
 $\qquad\qquad$ Surjeet Kour\footnote{Department of Mathematics, Indian Institute of Technology, Delhi, India. Email: koursurjeet@gmail.com.} \textrm{ }\footnote{Authors AB  and SK were partially supported by ANRF grant CRG/2023/004143 from the Government of India.} }

\date{}

\begin{document}

\maketitle

\begin{abstract}
We obtain Gerstenhaber type structures on  Davydov-Yetter cohomology with coefficients in half-braidings for a monoidal functor. Our approach uses a formal analogy between half-braidings of a monoidal functor and the entwining of a coalgebra with an algebra. We show that the Davydov-Yetter complex with coefficients carries the structure of a weak comp algebra. In particular, it is equipped with two distinct cup product structures $\cup$ and $\sqcup$ which are related in a manner that replaces graded commutativity. By considering entwining structures in the centralizer of the monoidal functor, we introduce subcomplexes of the  Davydov-Yetter complex with coefficients, which carry comp algebra structures. As a result, we obtain a graded commutative cup product and a graded Lie bracket on their cohomology which forms a Gerstenhaber algebra in the usual sense. 

\end{abstract}

\smallskip
{\bf MSC (2020) Subject Classification:} 16E40, 18M05.


\smallskip

\section{Introduction} 

The purpose of this paper is to study Gerstenhaber type structures, motivated by formal similarities between half-braidings in Davydov-Yetter cohomology with coefficients and the Hochschild cohomology theory of entwining structures. If $F:\mathscr C\longrightarrow \mathscr D$ is a monoidal functor between monoidal categories $\mathscr C$, $\mathscr D$, its deformations are described in terms of the Davydov-Yetter cohomology groups $H_{DY}^\bullet(F)$ (see, for instance, \cite{DB},\cite{DE}, \cite{Y1}, \cite{Y2}). We know that the Davydov-Yetter complex 
$C^\bullet_{DY}(F)$ of a monoidal functor $F$ is equipped with a number of interesting structures. In particular, $C^\bullet_{DY}(F)$ is equipped with a cup product
\begin{equation}
\cup : C^m_{DY}(F)\otimes C_{DY}^n(F)\longrightarrow C^{m+n}_{DY}(F)\qquad m,n\geq 0
\end{equation}
and  a graded  Lie bracket 
\begin{equation}
[\_\_,\_\_]: C^m_{DY}(F)\otimes C_{DY}^n(F)\longrightarrow C^{m+n-1}_{DY}(F)\qquad m,n\geq 0
\end{equation}
which together satisfy the conditions for $(C^\bullet_{DY}(F),\cup,[\_\_,\_\_])$ to be an $E_2$-algebra (see \cite[$\S$ 3]{DB}). At the level of cohomology, these operations descend to a Gerstenhaber algebra structure $(H^\bullet_{DY}(F),\cup,[\_\_,\_\_])$ on Davydov-Yetter cohomology of $F$. There is a well known analogy between the deformation cohomology of monoidal functors and the Hochschild cohomology of algebras. For an algebra $A$ over a field $k$, we know that the Hochschild cohomology groups $H^\bullet(A,A)$ of $A$ are equipped with a cup product $\cup$ and a Lie bracket $[\_\_,\_\_]$ which together determine a 
Gerstenhaber algebra $(H^\bullet(A,A),\cup, [\_\_,\_\_])$ (see, for instance, \cite[Theorem 1.4.9]{SJW}).

\smallskip
   In \cite{Gai}, Gainutdinov, Haferkamp and Schweigert  have introduced a framework for taking coefficients in the Davydov-Yetter cohomology of a monoidal functor $F:
   \mathscr C\longrightarrow \mathscr D$. Since $F$ is a monoidal functor, we have isomorphisms $\psi_{X,Y}: F(X)\otimes  F(Y)\overset{\cong}{\longrightarrow} F(X\otimes  Y)$ for
   $X$, $Y\in \mathscr C$. The coefficients for Davydov-Yetter cohomology of $F$ are taken from the collection $\mathcal Z(F)$ of half-braidings for the functor $F$ (see \cite[$\S$ 3]{Gai}): an object 
   $(U,\rho^U)\in \mathcal Z(F)$ consists of an object $U\in \mathscr D$ as well as a natural isomorphism   $\rho^U:U\otimes F(-)\overset{\cong}{\longrightarrow} F(-)\otimes U$ of functors  that is compatible with the isomorphisms $\psi_{X,Y}: F(X)\otimes F(Y)\overset{\cong}{\longrightarrow} F(X\otimes  Y)$. In this paper, we address the following question: to what extent do we have Gerstenhaber type structures on the Davydov-Yetter cohomology of a monoidal functor $F$ with coefficients in half-braidings $\mathcal Z(F)$ of $F$?
   
   \smallskip
We start in this paper by noting the formal similarity between half-braidings and the entwining of a coalgebra with an algebra. We recall (see \cite{BrMj}) that an entwining structure $(C,A,\phi)$ consists of a $k$-coalgebra $C$, a  $k$-algebra $A$ and a $k$-linear map $\phi:C\otimes A\longrightarrow A\otimes C$ satisfying certain conditions similar to the usual axioms for a bialgebra. For ease of comparison, we present side by side the commutative diagram appearing in  the half-braiding $(U,\rho^U)$ of a monoidal functor $F$ and that appearing in an entwining structure $(C,A,\phi)$. 
 \begin{equation}\label{12.2rc}
   \begin{array}{ccc} \xymatrix{
   U\otimes F(X)\otimes F(Y) \ar[rr]^{U\otimes \psi_{X,Y}} \ar[d]_{\rho^U(X)\otimes F(Y)} && U\otimes F(X\otimes Y) \ar[dd]^{\rho^U(X\otimes Y)}\\
   F(X)\otimes U\otimes F(Y) \ar[d]_{F(X)\otimes \rho^U(Y)}&&  \\
 F(X)\otimes F(Y)\otimes U  \ar[rr]^{\psi_{X,Y}\otimes U}&& F(X\otimes Y)\otimes U\\
   }&&  \xymatrix{
   C\otimes A \otimes A \ar[rr]^{C\otimes m} \ar[d]_{\phi\otimes A} && C\otimes A \ar[dd]^{\phi}\\
   A\otimes C\otimes A \ar[d]_{A\otimes \phi}&&  \\
A\otimes A\otimes C \ar[rr]^{m\otimes C}&\qquad \qquad & A\otimes C\\
   }\\
   \end{array}
   \end{equation} In \eqref{12.2rc}, $m:A\otimes A\longrightarrow A$ denotes the multiplication on the algebra $A$. 
Entwining structures were introduced by Brzezi\'{n}ski and Majid \cite{BrMj} as a kind of ``coalgebra principal bundle'' with applications to quantum homogenous spaces, coalgebra Galois extensions and gauge theory. In \cite{brz}, Brzezi\'{n}ski introduced a Hochschild cohomology theory for an entwining structure $(C,A,\phi)$, along with cup products  and other operations together satisfying a Gerstenhaber like structure known as a weak comp algebra. In this paper, we show that there are similar structures on the Davydov-Yetter cohomology of a monoidal functor $F$ taking certain kinds of coefficients in $\mathcal Z(F)$.  As such, we hope that this paper is the first step towards developing in detail the analogy between half-braidings and entwining structures. The literature on entwining structures is vast (see for instance, \cite{bbr}, \cite{bnkr}, \cite{brproc}, \cite{uni}, \cite{groot}, \cite{jia}, \cite{Sch}) and the modules over them, known as entwined modules (introduced in \cite{br3}), unify a number of classical concepts such as those of relative Hopf modules, Doi-Hopf modules, and Yetter-Drinfeld modules. Accordingly, one might expect in the future to also develop similar objects in the theory of monoidal functors, using half-braidings in place of entwining structures.

\smallskip
We now describe the paper in more detail. Given a monoidal functor $F:\mathscr C\longrightarrow \mathscr D$ between monoidal categories, we know (see \cite{Gai}, \cite{Majid}) that the collection $\mathcal Z(F)$ of half-braidings of $F$ carries the structure of a monoidal category. We fix a coalgebra object 
$(U,\rho^U,\Delta,\epsilon)$ as well as an algebra object $(R,\rho^R,\mu,\eta)$ in the monoidal category $\mathcal Z(F)$. Let $C^\bullet_{DY}(F,U,R)$ denote the Davydov-Yetter complex of $F$ with coefficients in $U$ and $R$, whose terms are given by (see \eqref{2.4e})
 \begin{equation}\label{2.4eIX}
   C^n_{DY}(F,U,R):=Nat(U\otimes F^{\otimes n},F^{\otimes n}\otimes R) \qquad n\geq 0
   \end{equation}  We note that the coefficients appearing in \eqref{2.4eIX} on the ``source side'' of the natural transformations are obtained from the  coalgebra object 
$(U,\rho^U,\Delta,\epsilon)$  in $\mathcal Z(F)$, while those on the ``target side'' are obtained from the algebra object  $(R,\rho^R,\mu,\eta)$ in $\mathcal Z(F)$. As with the Hochschild complex of an entwining structure developed in \cite{brz}, we show in Section 2 that the Davydov-Yetter complex 
$C^\bullet_{DY}(F,U,R)$ is equipped with not one but two distinct cup product structures
\begin{equation}\label{2cup2g}
\cup: C^m_{DY}(F,U,R)\otimes C^n_{DY}(F,U,R)\longrightarrow C^{m+n}_{DY}(F,U,R)\qquad \qquad \sqcup: C^m_{DY}(F,U,R)\otimes C^n_{DY}(F,U,R)\longrightarrow C^{m+n}_{DY}(F,U,R)\qquad m,n\geq 0
\end{equation} Further, the differential $\delta$ on the complex $C^\bullet_{DY}(F,U,R)$ is a graded derivation with respect to both the products in 
\eqref{2cup2g}, which leads to cup product structures $\cup$ and $\sqcup$ on the Davydov-Yetter cohomologies $H^\bullet_{DY}(F,U,R)$ of $F$ with
coefficients in $U$ and $R$ (see Proposition \ref{Pru2.4}). 

\smallskip
It is well known (see, for instance, \cite[Theorem 1.4.6]{SJW}) that the usual cup product on the Hochschild cohomology of an algebra is graded commutative. For Davydov-Yetter cohomology with coefficients, we will show that this graded commutativity is replaced by a certain relation between the cup products
$\cup$ and $\sqcup$ on 
$H^\bullet_{DY}(F,U,R)$. For this, we need the notion of a weak comp algebra introduced by  Brzezi\'{n}ski in \cite{brz}.  A weak comp algebra (see Definition \ref{wg}) consists of a graded vector space $V=\oplus_{m \geq 0} V^m$ and a collection of  $k$-linear maps 
\begin{equation} \diamond_i:V^m \otimes V^n \longrightarrow V^{m+n-1}\qquad \qquad m,n,i\geq 0
\end{equation} satisfying certain conditions obtained by weakening the axioms for a comp algebra. In particular, if  $V=\oplus_{m \geq 0} V^m$ is a comp algebra in the sense of Gerstenhaber and Schack \cite{GS92}, we know that $V^\bullet$ becomes a complex whose cohomology carries a cup product and a graded Lie bracket which together satisfy the conditions for being a Gerstenhaber algebra. In \cite{brz}, Brzezi\'{n}ski showed that the Hochschild complex of an entwining structure  forms a weak comp algebra. In Section 3, we introduce operations  (see \eqref{3.1q})
\begin{equation}\label{diam1c}
\diamond_i: C^m_{DY}(F,U,R) \otimes C^n_{DY}(F,U,R)\longrightarrow C^{m+n-1}_{DY}(F,U,R)\qquad m,n,i\geq 0
\end{equation} on the Davydov-Yetter complex $C^\bullet_{DY}(F,U,R)$. Then, the first main result of Section 3 is that $C^\bullet_{DY}(F,U,R)$ is a weak comp algebra (see Proposition \ref{p3.4f}).  Further, we show that the cup products $\cup$ and $\sqcup$, as well as the differential $\delta$ on $C^\bullet_{DY}(F,U,R)$ can be recovered in terms of the operations $\diamond_i$ in \eqref{diam1c}. The Davydov-Yetter complex $C^\bullet_{DY}(F,U,R)$ with coefficients in $U$ and $R$ then becomes a differential graded associative algebra with respect to both the products $\cup$ and $\sqcup$. On Davydov-Yetter cohomology groups $H^\bullet_{DY}(F,U,R)$ with coefficients,
the property of graded commutativity appearing in Hochschild cohomology  of an algebra is replaced by the following relation
\begin{equation} 
\bar{f} \cup \bar{g} = (-1)^{mn}~ \bar{g} \sqcup \bar{f} \qquad \bar{f} \in H^m_{DY}(F,U,R),\textrm{ } \bar{g} \in H^n_{DY}(F,U.R)
\end{equation} 
between the cup products $\cup$ and $\sqcup$ (see Theorem \ref{mdmsinhgr}). We give a number of examples of such situations, such as with group graded finite dimensional vector spaces. We can also consider Davydov-Yetter cohomology groups for an exact tensor functor $F:\mathscr C\longrightarrow \mathscr D$  between finite tensor categories, with coefficients in the Eilenberg-Moore category of modules over its central monad $Z_F$. In particular, we describe explicitly the cup products $\cup$ and $\sqcup$, as well as the operations $\diamond_i$ for the identity functor on the category $H-mod$ of finite dimensional modules over a finite dimensional Hopf algebra $H$.

\smallskip
In \cite{brz}, Brzezi\'{n}ski also introduced an equivariant subcomplex of the Hochschild complex of an entwining structure whose cohomology forms a Gerstenhaber algebra. In Section 4, we are inspired by this to study subcomplexes  
 of the Davydov-Yetter complex $C^\bullet_{DY}(F,U,R)$, which are ``equivariant'' with respect to the half-braidings on the coefficients, and whose cohomology can be equipped with the structure of a Gerstenhaber algebra. We do this in two steps. We first consider a  morphism $\,\overline{\!\Psi\!}\,: U \otimes R \longrightarrow R \otimes U$  in $\mathcal Z(F)$  which is well behaved with respect to the multiplication and the unit of the algebra object 
$R$. Corresponding to $\,\overline{\!\Psi\!}\,$, we construct subspaces $P^m(F,U,R,\,\overline{\!\Psi\!}\,)\subseteq C^m_{DY}(F,U,R)$, $m\geq 0$ consisting of natural transformations that satisfy a kind of equivariance condition (see \eqref{3.10cx}) with respect to the half-braiding $\rho^U$ on $U$. We show that $P^\bullet(F,U,R,\,\overline{\!\Psi\!}\,)$ forms a subcomplex of $C^\bullet_{DY}(F,U,R)$ and that the operations $\diamond_i$ on $C^\bullet_{DY}(F,U,R)$ restrict to $P^\bullet(F,U,R,\,\overline{\!\Psi\!}\,)$. Similarly, we can consider a morphism $\,\overline{\!\Phi\!}\,: U \otimes R \longrightarrow R \otimes U$ in  $\mathcal Z(F)$ which is well behaved with respect to  the comultiplication  and counit   of the coalgebra object $U$. Accordingly, we construct a subcomplex $Q^\bullet(F,U,R,\overline{\!\Phi\!}\,)\subseteq C^\bullet_{DY}(F,U,R)$ which consists of natural transformations  satisfying a kind of equivariance condition (see \eqref{4.5yy}) with respect to the half-braiding $\rho^R$ on $R$. 
 Further, the $\diamond_i$ operations on $C^\bullet_{DY}(F,U,R)$ restrict to $Q^\bullet(F,U,R,\overline{\!\Phi\!}\,)$.
 
 \smallskip
 We now put these two situations together and consider a morphism  $~\overline{\!\Theta\!}\,:U \otimes R \longrightarrow R \otimes U$ in $\mathcal{Z}(F)$ that is compatible with respect to both the coalgebra structure on $U$ and the algebra structure on $R$. In other words, $~\overline{\!\Theta\!}\,$ is an entwining structure in the monoidal category $\mathcal Z(F)$. Our main result in Section 4 is that  
 $P^\bullet(F,U,R,~\overline{\!\Theta\!}\,) $ and $Q^\bullet(F,U,R,~\overline{\!\Theta\!}\,)$ are weak comp algebras,  while the intersection
 \begin{equation}
D^\bullet(F,U,R,~\overline{\!\Theta\!}\,):= P^\bullet(F,U,R,~\overline{\!\Theta\!}\,) \cap Q^\bullet(F,U,R,~\overline{\!\Theta\!}\,)\subseteq C^\bullet_{DY}(F,U,R)
\end{equation}
 carries the structure of a comp algebra (see Theorem \ref{rysz}). 
  In particular, we show that $\cup$ restricts to a graded commutative cup product $\cup$ on the cohomology
groups  $H^\bullet(F, U, R, ~\overline{\!\Theta\!}\,)$  of the subcomplex $D^\bullet(F,U,R,~\overline{\!\Theta\!}\,)$.  There is also a graded Lie bracket
\begin{equation}\label{liein}
[f,g]:=f \diamond g - (-1)^{(m-1)(n-1)} g \diamond f\qquad \qquad f\in    D^m(F,U,R,~\overline{\!\Theta\!}\,), g\in D^n(F,U,R,~\overline{\!\Theta\!}\,)
\end{equation} where $f \diamond g:=\sum\limits_{i=0}^{m-1}  (-1)^{i(n-1)} f \diamond_i g$.   In Proposition \ref{P3.10x}, we show that the induced structures $\cup$ and $[\_\_,\_\_]$ together determine a  Gerstenhaber algebra structure on the cohomology  $H^\bullet(F, U, R, ~\overline{\!\Theta\!}\,)$. We also describe explicitly the complex $D^\bullet(F,U,R,~\overline{\!\Theta\!}\,)$ for certain entwining structures in the case of the identity functor on the category $H-mod$ of finite dimensional modules over a finite dimensional Hopf algebra $H$.

\smallskip
We also obtain a dual version of these constructions for an entwining structure in the opposite direction in the monoidal category $\mathcal Z(F)$. In other words, we consider a morphism 
$~\underline{\!\Theta\!}\,: R \otimes U \longrightarrow U \otimes R$  in $\mathcal Z(F)$ that is compatible with respect to the coalgebra structure on $U$ as well as the algebra structure on $R$. Accordingly, we obtain $P^\bullet(F,U,R,~\underline{\!\Theta\!}\,)\subseteq C^\bullet_{DY}(F,U,R)$ and $Q^\bullet(F,U,R,~\underline{\!\Theta\!}\,)\subseteq C^\bullet_{DY}(F,U,R)$ which are weak comp algebras, and we show that the intersection $D^\bullet(F,U,R,~\underline{\!\Theta\!}\,):= P^\bullet(F,U,R,~\underline{\!\Theta\!}\,) \cap Q^\bullet(F,U,R,~\underline{\!\Theta\!}\,)$ carries the structure of a comp algebra. We conclude by showing that the cohomology $H^\bullet(F,U,R,~\underline{\!\Theta\!}\,)$ of  $D^\bullet(F,U,R,~\underline{\!\Theta\!}\,)$ is equipped with a graded commutative cup product and a graded Lie bracket that gives $H^\bullet(F,U,R,~\underline{\!\Theta\!}\,)$ the structure of a Gerstenhaber algebra (see Theorem \ref{Th4.14ed}).  We would also like to thank the referee of an earlier version of this manuscript for making a number of insightful comments and suggestions.

   \section{Braidings of monoidal  functors and cup products}

   Let $k$ be a field and let $Vect_k$ be the category of finite dimensional vector spaces over $k$. By a monoidal category, we will always mean a $k$-linear abelian  strict monoidal category $(\mathscr C,\otimes_{\mathscr C},1_{\mathscr C})$ such that
   $\otimes_{\mathscr C}$ is  $k$-bilinear. Whenever the meaning is clear from context, we will drop the subscripts and write $(\mathscr C,\otimes_{\mathscr C},1_{\mathscr C})$ simply
   as $(\mathscr C,\otimes,1)$.  Let $\mathscr C$, $\mathscr D$ be monoidal categories and let $F:\mathscr C\longrightarrow \mathscr D$ be a monoidal  functor. This means that $F$ is a $k$-linear functor equipped with natural isomorphisms 
   \begin{equation}\label{002}
  \psi_{X,Y}: F(X)\otimes_{\mathscr D} F(Y)\overset{\cong}{\longrightarrow} F(X\otimes_{\mathscr C} Y)\qquad \eta: 1_{\mathscr D}\overset{\cong}{\longrightarrow} F(1_{\mathscr C})
   \end{equation} for $X$, $Y\in \mathscr C$ satisfying the usual associativity and unit conditions. 
   
   \begin{defn}\label{D2.1gt} (see  \cite[$\S$ 3]{Gai}) Let $F:\mathscr C\longrightarrow \mathscr D$ be a monoidal  functor. A half-braiding relative to $F$ consists of an object $U\in \mathscr D$ and a natural isomorphism 
   $\rho^U:U\otimes F(-)\overset{\cong}{\longrightarrow} F(-)\otimes U$ of functors such that the following diagram  
   \begin{equation}\label{2.2rc}
   \xymatrix{
   U\otimes F(X)\otimes F(Y) \ar[rr]^{U\otimes \psi_{X,Y}} \ar[d]_{\rho^U(X)\otimes F(Y)} && U\otimes F(X\otimes Y) \ar[dd]^{\rho^U(X\otimes Y)}\\
   F(X)\otimes U\otimes F(Y) \ar[d]_{F(X)\otimes \rho^U(Y)}&&  \\
 F(X)\otimes F(Y)\otimes U  \ar[rr]^{\psi_{X,Y}\otimes U}&& F(X\otimes Y)\otimes U\\
   }
   \end{equation}
   is commutative for every $X$, $Y\in \mathscr C$. The collection of such pairs $(U,\rho^U)$ forms a category,  which is known as the centralizer $\mathcal Z(F)$ of the functor $F:
   \mathscr C\longrightarrow \mathscr D$. A morphism 
   $f:(U,\rho^U)\longrightarrow (V,\rho^V)$ in the category $\mathcal Z(F)$ consists of a morphism $f:U\longrightarrow V$ in $\mathscr D$ such that $(F(X)\otimes f)~\circ \rho^U(X)=\rho^V(X)~\circ 
   (f\otimes  F(X)):U\otimes F(X)\longrightarrow F(X)\otimes V$ for each $X\in \mathscr C$
   \end{defn} 
   
  We know that  the monoidal structure on $\mathscr D$ induces a monoidal structure on the category 
   $\mathcal Z(F)$ (see \cite{Gai}, \cite{Majid}) as follows: if $(U,\rho^U)$, $(V,\rho^V)\in \mathcal Z(F)$, we have
   \begin{equation}\label{monz1}
   (U,\rho^U)\otimes_{\mathcal Z(F)}(V,\rho^V):=\left(U\otimes V,\rho^{U\otimes V}=\left\{\rho^{U\otimes V}(X):U\otimes V\otimes F(X)\underset{\cong}{\xrightarrow{U\otimes \rho^V(X)}}U\otimes F(X)\otimes V
   \underset{\cong}{\xrightarrow{\rho^U(X)\otimes V}}F(X)\otimes U\otimes V\right\}_{X\in Ob(\mathscr C)}\right)
   \end{equation} The unit object in $\mathcal Z(F)$ is given by the pair $\left( 1_{\mathscr D}, 1_{\mathscr D}\otimes  F(-)\overset{\cong}{\longrightarrow} F(-)\otimes 1_{\mathscr D}\right)$.
   
   \smallskip
 From now onwards, we will assume that 
  the monoidal  functor $F:\mathscr C\longrightarrow \mathscr D$ preserves direct sums (whenever they exist).   For $n\geq 0$, let 
   $F^{\otimes n}$ denote the functor
   \begin{equation}
   F^{\otimes n}: \underbrace{\mathscr C \times ... \times \mathscr C}_{\mbox{\small $n$-times}} \longrightarrow \mathscr D  \qquad (X_1,...,X_n)\mapsto F(X_1)\otimes ...\otimes F(X_n) 
   \end{equation} where $F^{\otimes 0}$ is understood to be the functor $F^{\otimes 0}:Vect_k\longrightarrow \mathscr D$ that preserves direct sums and takes the vector space 
   $k$ to $1_{\mathscr D}\in \mathscr D$. 
   
   \smallskip
   In \cite{Gai}, Gainutdinov, Haferkamp and Schweigert have shown that the collection $\mathcal Z(F)$ of half-braidings gives a framework for introducing coefficients in the Davydov-Yetter cohomology  of the
   functor $F$. Accordingly (see \cite[$\S$ 4.2]{Fai}, \cite[$\S$ 3]{Gai}), the Davydov-Yetter cohomology of $F$ with coefficients in
   half-braidings $(U,\rho^U)$ and $(R,\rho^R)$ is obtained from a complex $C^\bullet_{DY}(F,U,R)$ whose terms are
   \begin{equation}\label{2.4e}
   C^n_{DY}(F,U,R):=Nat(U\otimes F^{\otimes n},F^{\otimes n}\otimes R) \qquad n\geq 0
   \end{equation} 
   By definition, a natural transformation $f\in  C^n_{DY}(F,U,R)$ can be described by a family of morphisms
   \begin{equation}\label{2.6e}
  f=\{ f_{X_1,...,X_n}: U\otimes F(X_1)\otimes ... \otimes F(X_n)\longrightarrow F(X_1)\otimes ...\otimes F(X_n) \otimes R  \}_{(X_1,...,X_n)\in Ob(\mathscr C)^n}
   \end{equation}  For ease of notation, we always assume that  $F:\mathscr{C} \longrightarrow \mathscr{D}$  is strict monoidal, i.e.,  $\psi_{X,Y}=\mathrm{id}$ for all
   $X$, $Y\in Ob(\mathscr C)$. However, the definitions and results can be extended to non-strict functors as well by properly inserting $\psi$. 
   
\smallskip 
  For $n>0$, the differential $\partial^n:C^n_{DY}(F,U,R)\longrightarrow C^{n+1}_{DY}(F,U,R)$ is given by $\partial^n=\sum_{i=0}^{n+1} (-1)^i \partial_i^n$, where the $\partial_i^n$ are defined as follows:
  
\begin{equation}\label{dif}
\begin{aligned}
\partial_0^n(f)_{X_1,\ldots,X_{n+1}}
&=
\bigl({F(X_1)} \otimes f_{X_2,\ldots,X_{n+1}}\bigr)
\circ
\bigl(\rho^{U}(X_1) \otimes {F(X_2) \otimes \ldots \otimes F(X_{n+1})}\bigr),
\\
\partial_i^n(f)_{X_1,\ldots,X_{n+1}}
&=
f_{X_1,\ldots,X_{i} \otimes X_{i+1},\ldots,X_{n+1}}
\qquad \text{for } 1 \le i \le n,\\
\partial_{n+1}^n(f)_{X_1,\ldots,X_{n+1}}
&=
\bigl( {F(X_1) \otimes \cdots \otimes F(X_{n})} \otimes \rho^{R}({X_{n+1}})\bigr)
\circ
\bigl(f_{X_1,\ldots,X_{n}} \otimes  {F(X_{n+1})}\bigr).
\end{aligned}
\end{equation}

 For $n=0$, we have $\partial^0(f)_{X_1}:=(F(X_1)\otimes f)\circ \rho^U(X_1)- \rho^R(X_1) \circ (f\otimes F(X_1))$. The cohomology groups of the complex $C^\bullet_{DY}(F,U,R)$ will be denoted by $H^\bullet_{DY}(F,U,R)$. 
   
   \smallskip  For any ordered tuple $(X_1, \ldots, X_n)$ of objects in $\mathscr C$, we set $\rho^U(X_1, \ldots, X_n)$ to be the composition
  \begin{equation*}\small
  \xymatrix{
  U\otimes F(X_1)\otimes \ldots \otimes F(X_n)\ar[d]_{\rho^U(X_1)\otimes F(X_2)\otimes \ldots \otimes F(X_n)}\ar[rrrrrrrr]^{ \rho^U(X_1, \ldots, X_n)} &&&&&&&&  F(X_1)\otimes \ldots \otimes F(X_n)\otimes U\\
  F(X_1)\otimes U\otimes F(X_2)\otimes \ldots \otimes F(X_n)\ar[rrrr]^{\qquad \qquad \qquad F(X_1)\otimes \rho^U(X_2)\otimes F(X_3)\otimes \ldots\otimes F(X_n)} &&&& \ldots \ar[rrrr]^{F(X_1)\otimes \ldots
  \otimes F(X_{n-2})\otimes \rho^U(X_{n-1})\otimes F(X_n)\qquad \qquad \qquad   } &&&& \ar[u]_{F(X_1)\otimes \ldots \otimes F(X_{n-1})\otimes \rho^U(X_n)} F(X_1)\otimes \ldots \otimes F(X_{n-1})\otimes U\otimes F(X_n)
  }
  \end{equation*} 
  
  From now onwards, we suppose that $(U,\rho^U)$ is additionally equipped with the structure  $(U,\rho^U, \Delta,\epsilon)$ of a coalgebra object in the monoidal category
  $\mathcal Z(F)$, while $(R,\rho^R)$ is equipped with the structure $(R,\rho^R,\mu,\eta)$ of an algebra object.  Since $\Delta$, $\epsilon$, $\mu$ and $\eta$ are morphisms in $\mathcal Z(F)$, it follows from the 
    monoidal structure on 
$\mathcal Z(F)$ described in \eqref{monz1} that the following diagrams commute for each $X\in \mathscr C$
   \begin{equation}\label{2.8bv}
   \begin{array}{cc}
   \xymatrix{
   U\otimes F(X)\ar[d]_{\Delta\otimes F(X)} \ar[rrrr]^{\rho^U(X)} &&&& F(X)\otimes U\ar[d]^{F(X)\otimes\Delta } \\
U\otimes U\otimes F(X) \ar[rr]^{U\otimes {\rho^U(X)}}  &&U\otimes F(X)\otimes U\ar[rr]^{{\rho^U(X)}\otimes U}&& F(X)\otimes U\otimes U\\
   }
   &
   \xymatrix{
   U\otimes F(X) \ar[rr]^{\rho^U(X)}\ar[drr]_{\epsilon\otimes F(X)}&& F(X)\otimes U \ar[d]^{F(X)\otimes \epsilon}\\
   && F(X)\\
   }\\
   \end{array}
   \end{equation} 
   \begin{equation}\label{Rcom}
   \begin{array}{cc}
   \xymatrix{
    R \otimes R \otimes F(X) \ar[d]_{\mu \otimes F(X)} \ar[rr]^{R \otimes \rho^R(X) } && R \otimes F(X) \otimes R \ar[rr]^{\rho^R(X) \otimes R} && F(X) \otimes R \otimes R  \ar[d]^{{F(X)} \otimes \mu}\\
  R \otimes F(X)  \ar[rrrr]^{\rho^R(X)} &&&& F(X) \otimes R
   }
   &
   \xymatrix{
   F(X) \ar[d]_{\eta\otimes F(X)} \ar[drr]^{F(X)\otimes \eta}\\
   R\otimes F(X)   \ar[rr]^{\rho^R(X)} &&   F(X)\otimes R
   }
   \end{array}
   \end{equation} 
  
 \smallskip 
  We are now ready to define the first cup product on $C^\bullet_{DY}(F,U,R)$.  For $m$, $n\geq 0$, we set
   \begin{equation}\label{cup2}
   \begin{array}{c}
   \cup: C^m_{DY}(F,U,R)\otimes C^n_{DY}(F,U,R)\longrightarrow C^{m+n}_{DY}(F,U,R)\qquad f\otimes g\mapsto f\cup g
   \end{array}
   \end{equation}
with $(f \cup g)_{X_1,...,X_{m+n}}$ to be the following composition

\smallskip
 \begin{equation}\small
 (f \cup g)_{X_1,...,X_{m+n}}:=\left\lbrace \begin{array}{cc}
 \xymatrix{
 U \otimes F(X_1) \otimes \ldots \otimes F(X_m) \otimes F(X_{m+1}) \otimes \ldots \otimes F(X_{m+n}) \ar[d]^{\Delta \otimes F(X_1) \otimes \ldots \otimes F(X_{m+n})} \\
 U \otimes U \otimes F(X_1) \otimes \ldots \otimes F(X_m) \otimes F(X_{m+1}) \otimes \ldots \otimes F(X_{m+n})
 \ar[d]^{U \otimes \rho^U(X_1,\ldots,X_m) \otimes F(X_{m+1}) \otimes \ldots \otimes F(X_{m+n})} \\
 U \otimes F(X_1) \otimes \ldots \otimes F(X_m) \otimes U \otimes F(X_{m+1}) \otimes \ldots \otimes F(X_{m+n}) \ar[d]^{f_{X_1,\ldots,X_m} \otimes g_{X_{m+1},\ldots,X_{m+n}}} \\
 F(X_1) \otimes \ldots \otimes F(X_m) \otimes R \otimes F(X_{m+1}) \otimes \ldots \otimes F(X_{m+n}) \otimes R  \ar[d]^{F(X_1) \otimes \ldots \otimes F(X_m) \otimes \rho^R(X_{m+1},\ldots,X_{m+n}) \otimes R} \\
 F(X_1) \otimes \ldots \otimes F(X_m) \otimes F(X_{m+1}) \otimes \ldots \otimes F(X_{m+n}) \otimes R \otimes R \ar[d]^{F(X_1) \otimes \ldots F(X_{m+n}) \otimes \mu}\\
 F(X_1) \otimes \ldots \otimes F(X_{m+n}) \otimes R 
 }
 \end{array}\right\rbrace
 \end{equation}
 

for any $X_1, \ldots, X_{m+n}\in Ob(\mathscr C)$.  
We first show that the cup product in \eqref{cup2} induces a product on the cohomology $C^\bullet_{DY}(F,U,R)$. 

   \begin{lem}\label{L2.1} For $f\in C^m_{DY}(F,U,R)$ and $g\in C^n_{DY}(F,U,R)$, we have
   \begin{equation}\label{2.10u}
   \partial^{m+n}(f\cup g)=(\partial^m(f)\cup g)+(-1)^m(f\cup \partial^n(g))
   \end{equation}
   \end{lem}
   \begin{proof}  Using the definition of differential in \eqref{dif}, we first observe that for $m,n >0$

\begin{equation}\label{2.13dx}
\partial_i^{m+n}(f \cup g) =
\begin{cases}
\partial_i^m(f) \cup g, & \text{if } 0 \le i \le m, \\[4pt]
f \cup \partial_{i-m}^n(g), & \text{if } m < i \le m+n+1
\end{cases}
\end{equation}

and also that \begin{equation}\label{2.14dx}
\partial_{m+1}^m(f) \cup g = f \cup \partial_0^n(g)
\end{equation}
We now see that 

   \begin{align*}
\partial^{m+n}(f \cup g)
&~~= \sum_{i=0}^{m+n+1} (-1)^i \partial_i^{m+n}(f \cup g) \stackrel{\eqref{2.13dx}}{=}
\left(
\sum_{i=0}^{m} (-1)^i \partial_i^m(f) \cup g
\right)
+
\left(
\sum_{i=m+1}^{m+n+1} (-1)^i f \cup \partial_{i-m}^n(g)
\right) \\
& \stackrel{\eqref{2.14dx}}{=}
\left(
\sum_{i=0}^{m} (-1)^i \partial_i^m(f) \cup g
\right)
+ (-1)^{m+1} \partial_{m+1}^m(f) \cup g
+ (-1)^m f \cup \partial_0^n(g) + \left(
\sum_{i=1}^{n+1} (-1)^{m+i} f \cup \partial_i^n(g)
\right) \\
& ~~=
\left(
\sum_{i=0}^{m+1} (-1)^i \partial_i^m(f) \cup g
\right)
+
\left(
\sum_{i=0}^{n+1} (-1)^{m+i} f \cup \partial_i^n(g)
\right) = \partial^m(f) \cup g + (-1)^m f \cup \partial^n(g).
\end{align*} The case where $m=0$ or $n=0$ may be verified directly, by making use of the commutative diagrams \eqref{2.8bv} and \eqref{Rcom}. From this, the result follows.

\end{proof}

We  now proceed to define the second cup product on $C^\bullet_{DY}(F,U,R)$. For $m$, $n\geq 0$, we set 
 \begin{equation}\label{cup3}  
   \sqcup: C^m_{DY}(F,U,R)\otimes C^n_{DY}(F,U,R)\longrightarrow C^{m+n}_{DY}(F,U,R)\qquad f\otimes g\mapsto f\sqcup g\\ 
\end{equation}   with $  (f\sqcup g)_{X_1,...,X_{m+n}}$ to be the following composition

\begin{equation}
\small
(f\sqcup g)_{X_1, \ldots, X_{m+n}}:=\left\lbrace \begin{CD}
U\otimes F(X_1)\otimes \ldots \otimes F(X_m) \otimes F(X_{m+1})\otimes \ldots F(X_{m+n})\\
@VV\Delta\otimes F(X_1)\otimes \ldots \otimes F(X_{m+n})V \\
U\otimes U\otimes F(X_1)\otimes \ldots \otimes F(X_m) \otimes F(X_{m+1}) \otimes \ldots \otimes F(X_{m+n})\\
@VVU\otimes f_{X_1,\ldots, X_{m}}\otimes F(X_{m+1})\otimes \ldots \otimes F(X_{m+n})V \\
U\otimes F(X_1)\otimes \ldots \otimes F(X_m) \otimes R \otimes F(X_{m+1}) \otimes \ldots \otimes F(X_{m+n})\\
@VV\rho^U(X_1, \ldots, X_m)\otimes \rho^R(X_{m+1}, \ldots,  X_{m+n})V \\
F(X_{1})\otimes \ldots \otimes F(X_{m})\otimes U \otimes F(X_{m+1})\otimes \ldots \otimes F(X_{m+n}) \otimes R\\
@VVF(X_{1})\otimes \ldots \otimes F(X_{m})\otimes g_{X_{m+1}, \ldots, X_{m+n}} \otimes RV\\
F(X_{1})\otimes \ldots \otimes F(X_m) \otimes F(X_{m+1}) \otimes \ldots \otimes F(X_{m+n}) \otimes R \otimes R \\
@VV {F(X_1) \otimes \ldots  \otimes F(X_{m+n})} \otimes \mu V\\
F(X_{1})\otimes \ldots \otimes F(X_m) \otimes F(X_{m+1})\otimes \ldots \otimes F(X_{m+n}) \otimes R
\end{CD} \right\rbrace
\end{equation} for $X_1, \ldots, X_{m+n}\in Ob(\mathscr C)$.

   \begin{lem}\label{L2.2} For $f\in C^m_{DY}(F,U,R)$, $g\in C^n_{DY}(F,U,R)$, we have
   \begin{equation}\label{2.14c}
   \partial^{m+n}(f\sqcup g)=\partial^m(f)\sqcup g+(-1)^m(f\sqcup \partial^n(g))
   \end{equation}
   \end{lem}
   \begin{proof}
The proof  follows in a manner similar to Lemma \ref{L2.1} by using \eqref{2.8bv}, \eqref{Rcom} and by observing that
\[
\partial_i^{m+n}(f \sqcup g) =
\begin{cases}
\partial_i^m(f) \sqcup g, & \text{if } 0 \le i \le m, \\[4pt]
f \sqcup \partial_{i-m}^n(g), & \text{if } m < i \le m+n+1
\end{cases}
\]

as well as that
$\partial_{m+1}^m(f) \sqcup g = f \sqcup \partial_0^n(g)$.

\end{proof}

\begin{prop}\label{Pru2.4}
Let $F:\mathscr C\longrightarrow \mathscr D$ be a  monoidal  functor between monoidal  categories. Let $(U,\rho^U,\Delta,\epsilon)$ be a coalgebra object and $(R, \rho^R, \mu,\eta)$ be an algebra object in the centralizer $\mathcal Z(F)$ of $F$. Then, the cup products $\cup$ as in \eqref{cup2} and $\sqcup$ as in \eqref{cup3} induce products on the Davydov-Yetter cohomology $H^{\bullet}_{DY}(F,U,R)$.
\begin{equation*}
\cup : H^m_{DY}(F,U,R)\otimes H^n_{DY}(F,U,R)\longrightarrow H^{m+n}_{DY}(F,U,R) \qquad \qquad
\sqcup: H^m_{DY}(F,U,R)\otimes H^n_{DY}(F,U,R)\longrightarrow H^{m+n}_{DY}(F,U,R) 
\end{equation*}
for $m$, $n\geq 0$. 
\end{prop}
\begin{proof}
This follows from  Lemma \ref{L2.1} and Lemma \ref{L2.2}.
\end{proof}	

\section{Weak comp algebras and the  Davydov-Yetter complex}

We continue with a (strict) monoidal  functor  $F:\mathscr C\longrightarrow \mathscr D$, a coalgebra object $(U,\rho^U,\Delta,\epsilon)$ and an algebra object $(R,\rho^R, \mu,\eta)$  in the centralizer $\mathcal Z(F)$ of $F$.
In this section, we will  study Gerstenhaber algebra like structures on the complex $C^\bullet_{DY}(F,U,R)$ and its cohomology. We begin by recalling the notion of a weak comp algebra due to Brzezi{\'n}ski \cite[Definition 4.4]{brz}. 

\begin{defn}\label{wg}
Let $k$ be a field. A (right) weak comp algebra over $k$ consists of the following data: (a) a graded $k$-vector space $V=\oplus_{m \geq 0} V^m$,
(b) an element $\pi \in V^2$, and (c) a collection of  $k$-linear maps 
\begin{equation} \diamond_i:V^m \otimes V^n \longrightarrow V^{m+n-1} \qquad i\geq 0
\end{equation} such that the following conditions are satisfied for any $f \in V^m, g \in V^n$ and $h \in V^p$

\smallskip
(1) $f \diamond_i g=0$ for $i > m-1$

\smallskip
(2) $(f \diamond_i g) \diamond_j h=f \diamond_i (g \diamond_{j-i} h)$ if $i \leq j < n+i$

\smallskip
(3) if either $g=\pi$ or $h=\pi$,

\vspace{.1cm}
~~~~~~~~$(f \diamond_i g) \diamond_j h=(f \diamond_j h) \diamond_{i+p-1} g$ if $j < i$

\smallskip
(4) $\pi \diamond_0 \pi=\pi \diamond_1 \pi$
\end{defn}

A weak comp algebra is a comp algebra in the sense of Gerstenhaber and Schack \cite{GS92} if condition (3) in Definition \ref{wg} holds for all $g\in V^n$ and $h\in V^p$.  For more on such structures and their role in the Hochschild cohomology of algebras, we refer the reader, for instance, to \cite{cgs}, \cite{Ger}.

\smallskip
Our first aim in this section is to show that the complex $C^\bullet_{DY}(F,U,R)$ is a weak comp algebra. For this, the  main step  is to define operations
\begin{equation}
\diamond_i: C^m_{DY}(F,U,R) \otimes C^n_{DY}(F,U,R)\longrightarrow C^{m+n-1}_{DY}(F,U,R)\qquad i\geq 0
\end{equation} For $f \in C^m_{DY}(F,U,R)=Nat(U \otimes F^{\otimes m}, F^{\otimes m} \otimes R)$, $g \in C^n_{DY}(F,U,R)=Nat(U \otimes F^{\otimes n}, F^{\otimes n} \otimes R)$ and $0 \leq i <m$, we set $f \diamond_i g \in C^{m+n-1}_{DY}(F,U,R)=Nat(U \otimes F^{\otimes m+n-1}, F^{\otimes m+n-1} \otimes R)$ to be the natural transformation determined by the family $
f \diamond_i g :=\{ (f \diamond_i g)_{X_1 \ldots, X_{m+n-1}} \}_{X_j  \in Ob(\mathscr C), 1 \leq j \leq m+n-1}
$
where each   $(f \diamond_i g)_{X_1 \ldots, X_{m+n-1}}$ is set to be the following composition

\begin{equation}\label{3.1q}\small 
(f \diamond_i g)_{X_1 \ldots, X_{m+n-1}}:=\left\lbrace  \begin{tikzcd}
U\otimes F(X_1) \otimes \ldots \otimes F(X_{m+n-1})
\arrow[d,"{\Delta \otimes F(X_1) \otimes \ldots \otimes F(X_{m+n-1}) }"]\\
U \otimes U \otimes F(X_1) \otimes \ldots \otimes F(X_{m+n-1})
\arrow[d,"{U \otimes \rho^U(X_1, \ldots, X_i) \otimes F(X_{i+1}) \otimes \ldots \otimes  F(X_{m+n-1}) }"]\\
U \otimes F(X_1) \otimes \ldots \otimes F(X_i) \otimes U \otimes F(X_{i+1}) \otimes \ldots \otimes  F(X_{m+n-1})
\arrow[d,"{ U \otimes F(X_1) \otimes \ldots \otimes F(X_i) \otimes g_{X_{i+1}, \ldots ,X_{i+n}} \otimes F(X_{i+n+1})\otimes \ldots\otimes  F(X_{m+n-1}) }"]\\
U \otimes F(X_1) \otimes \ldots \otimes F(X_i) \otimes F(X_{i+1} \otimes \ldots \otimes X_{i+n}) \otimes R \otimes F(X_{i+n+1}) \otimes \ldots \otimes  F(X_{m+n-1})
\arrow[d,"{U \otimes F(X_1) \otimes \ldots \otimes F(X_i) \otimes F(X_{i+1} \otimes \ldots \otimes X_{i+n}) \otimes \rho^R(X_{i+n+1},\ldots, X_{m+n-1})  }"]\\
U \otimes F(X_1) \otimes \ldots \otimes F(X_i) \otimes F(X_{i+1} \otimes \ldots \otimes X_{i+n}) \otimes F(X_{i+n+1}) \otimes \ldots \otimes F(X_{m+n-1}) \otimes R
\arrow[d,"{f_{X_1, \ldots ,X_i, X_{i+1} \otimes \ldots \otimes X_{i+n}, X_{i+n+1}, \ldots, X_{m+n-1}} \otimes R }"]\\
F(X_1) \otimes \ldots \otimes F(X_i) \otimes F(X_{i+1} \otimes \ldots \otimes X_{i+n}) \otimes F(X_{i+n+1}) \otimes \ldots \otimes F(X_{m+n-1}) \otimes R \otimes R
\arrow[d,"{F(X_1) \otimes \ldots \otimes F(X_{m+n-1}) \otimes \mu}"]\\
F(X_1) \otimes \ldots \otimes F(X_{m+n-1}) \otimes R 
\end{tikzcd}\right\rbrace
\end{equation}

For all other values of $i$, we set $f \diamond_i g=0$. We also define $\pi:U \otimes F(-) \otimes F(-) \longrightarrow F(-) \otimes F (-)\otimes R\in C^2_{DY}(F,U,R)$ as
\begin{align}\label{pii}
\pi:=\{\pi_{X_1,X_2}: \xymatrixcolsep{5pc}\xymatrix{U \otimes F(X_1) \otimes F(X_2) \ar[r]^{\epsilon \otimes F(X_1) \otimes F(X_2)} &  F(X_1) \otimes F(X_2) \ar[r]^{F(X_1) \otimes F(X_2) \otimes \eta} &  F(X_1) \otimes F(X_2) \otimes R}\}_{X_1,X_2 \in \text{Ob}(\mathscr C)}
\end{align} We now have the following result.

\begin{lem}\label{con2} Let $f \in  C^m_{DY}(F,U,R)$, $g \in  C^n_{DY}(F,U,R)$ and $h \in  C^p_{DY}(F,U,R)$. Then,  
\begin{itemize}
\item[(a)]  for $i \leq j < n+i$, we have
\begin{equation}\label{3.4tt}
(f \diamond_i g) \diamond_j h=f \diamond_i (g \diamond_{j-i} h) \in  C^{m+n+p-2}_{DY}(F,U,R)
\end{equation}
\item[(b)] if $g = \pi$ or $h=\pi$ and $j < i$, we have
 $
(f \diamond_i g) \diamond_j h=(f \diamond_j h) \diamond_{i+p-1} g 
$. 
\end{itemize}
\end{lem}
\begin{proof}
We only show in detail the computation for  (a). Part  (b) may be verified in a similar manner. 
It is clear that \eqref{3.4tt} is satisfied for all $i \geq m$. Hence, we assume $0 \leq i < m$. 
Let $X_1, \ldots, X_{m+n+p-2}\in Ob(\mathscr C)$. Then,
using the compatibility between $\rho^U$ and $\Delta$ (resp. $\rho^R$ and $\mu$) as in \eqref{2.8bv} (resp. \eqref{Rcom}), the right hand side of \eqref{3.4tt} fits into the following commutative diagram:

\begin{center}
\scalebox{0.6}{
\begin{tikzcd}[ampersand replacement=\&, column sep=2.5em, row sep=3em, nodes={inner sep=2pt}] 
\boxed{\begin{tabular}{l} $U \otimes F(X_1) \otimes \ldots$ \\ $\ldots \otimes F(X_{m+n+p-2})$ \end{tabular}} 
\arrow[rr, "{\Delta \otimes \ldots}"] 
\arrow[d, "{\Delta \otimes \ldots}"] 
\& \& \boxed{\begin{tabular}{l} $U \otimes U \otimes F(X_1) \otimes \ldots$ \\ $\ldots \otimes F(X_{m+n+p-2})$ \end{tabular}} 
\arrow[rr, "{\_ \otimes \rho^U(X_1,\ldots,X_i) \otimes \ldots}"] 
\& \& \boxed{\begin{tabular}{l} $U \otimes F(X_1) \otimes \ldots \otimes F(X_i) \otimes U \otimes$\\ $F(X_{i+1}) \otimes \ldots \otimes F(X_{m+n+p-2})$ \end{tabular}} 
\arrow[d, xshift=0pt, "{\_ \otimes \ldots \otimes \Delta \otimes \ldots}"] \\ 
\boxed{\begin{tabular}{l} $U \otimes U \otimes F(X_1) \otimes \ldots$ \\ $\ldots \otimes F(X_{m+n+p-2})$ \end{tabular}} 
\arrow[d,  yshift=3pt, "{\_ \otimes \Delta \otimes \ldots}"]
\& \& \& \& \boxed{\begin{tabular}{l} $U \otimes F(X_1) \otimes \ldots \otimes F(X_i)\otimes U \otimes U \otimes$\\ $F(X_{i+1}) \ldots \otimes F(X_{m+n+p-2})$ \end{tabular}} 
\arrow[d, "{\_ \otimes \_ \otimes \ldots \otimes \_ \otimes \rho^U(X_{i+1},\ldots,X_j) \otimes \ldots}"] \\
\boxed{\begin{tabular}{l} $U \otimes U \otimes U \otimes F(X_1) \otimes \ldots$ \\ $\ldots \otimes F(X_{m+n+p-2})$ \end{tabular} 
\arrow[d,"{\_  \otimes \rho^U(X_1,\ldots,X_i) \otimes \ldots}"] } 
\& \& \& \& \boxed{\begin{tabular}{l} $U \otimes F(X_1) \otimes \ldots \otimes F(X_i) \otimes U \otimes F(X_{i+1}) \otimes $\\ $\ldots \otimes F(X_j) \otimes U \otimes \ldots \otimes F(X_{m+n+p-2})$ \end{tabular}} 
\arrow[d, "{\_ \otimes \ldots \otimes \_ \otimes \ldots \otimes h_{X_{j+1}, \ldots ,X_{j+p}} \otimes \ldots}"] \\ 
\boxed{\begin{tabular}{l} $U \otimes U \otimes F(X_1) \otimes \ldots \otimes F(X_i) \otimes U$\\ $\otimes \ldots \otimes F(X_{m+n+p-2})$ \end{tabular} } 
\arrow[d,"{ \_ \otimes \rho^U(X_1,\ldots,X_i) \otimes \_ \otimes \ldots }"] 
\& \& \& \& 
\boxed{\begin{tabular}{l} $U \otimes F(X_1) \otimes \ldots \otimes F(X_i) \otimes U \otimes F(X_{i+1}) \otimes \ldots \otimes F(X_j) \otimes $\\ $F(X_{j+1}) \otimes \ldots \otimes F(X_{j+p}) \otimes R \otimes \ldots \otimes F(X_{m+n+p-2})$ \end{tabular}} 
\arrow[d, "{\_ \otimes \ldots \otimes \_ \otimes \ldots \otimes \rho^R(X_{j+p+1}, \ldots, X_{n+p+i-1}) \otimes \ldots}"]\\ 
\boxed{\begin{tabular}{l} $U \otimes F(X_1) \otimes \ldots \otimes F(X_i) \otimes U \otimes U$ \\ 
$\otimes F(X_{i+1}) \otimes \ldots \otimes F(X_{m+n+p-2})$ \end{tabular} } 
\arrow[d,"{\_ \otimes \ldots \otimes \_ \otimes \rho^U(X_{i+1},\ldots,X_j) \otimes \ldots }"] 
\& \& \& \& 
\boxed{\begin{tabular}{l} $U \otimes F(X_1) \otimes...\otimes F(X_i) \otimes U \otimes F(X_{i+1}) \otimes... \otimes F(X_{j}) \otimes$\\ $ F(X_{j+1} \otimes \ldots \otimes X_{j+p}) \otimes F(X_{j+p+1}) \otimes \ldots \otimes F(X_{n+p+i-1}) \otimes $\\ $R \otimes \ldots \otimes F(X_{m+n+p-2}) $ \end{tabular}} 
\arrow[d,"{\_ \otimes \ldots \otimes g_{X_{i+1},\ldots,X_j,X_{j+1} \otimes \ldots X_{j+p},\ldots X_{n+p+i-1})} \otimes \_ \otimes \ldots }"]\\ 
\boxed{\begin{tabular}{l} $U \otimes F(X_1) \otimes \ldots \otimes F(X_i) \otimes U \otimes $ \\ $F(X_{i+1}) \otimes \ldots F(X_j) \otimes U\otimes \ldots \otimes F(X_{m+n+p-2})$ \end{tabular} } 
\arrow[d,"{\_ \otimes \ldots  \otimes \_ \otimes \ldots \otimes h_{X_{j+1},\ldots,X_{j+p}}} \otimes \ldots"] 
\& \& \& \& 
\boxed{\begin{tabular}{l} $U \otimes F(X_1) \otimes...\otimes F(X_i) \otimes F(X_{i+1}) \otimes... \otimes F(X_{j}) \otimes$\\ $ F(X_{j+1} \otimes \ldots \otimes X_{j+p}) \otimes F(X_{j+p+1}) \otimes \ldots \otimes F(X_{n+p+i-1}) \otimes R \otimes R$\\ $\otimes \ldots \otimes F(X_{m+n+p-2}) $ \end{tabular}} 
\arrow[d,"{\_ \otimes \ldots \otimes \mu \otimes \ldots}"]\\ 
\boxed{ \begin{tabular}{l} $U \otimes F(X_1) \otimes...\otimes F(X_i) \otimes U \otimes F(X_{i+1}) \otimes \ldots$\\ $\otimes F(X_{j+1} \otimes \ldots \otimes X_{j+p}) \otimes R \otimes \ldots \otimes F(X_{m+n+p-2})$ \end{tabular} } 
\arrow[d,"{\_ \otimes \ldots \_ \otimes \ldots \otimes \rho^R(X_{j+p+1}, \ldots, X_{n+p+i-1}) \otimes \ldots }"] 
\& \& \& \& 
\boxed{\begin{tabular}{l} $U \otimes F(X_1) \otimes...\otimes F(X_i) \otimes F(X_{i+1}) \otimes... \otimes F(X_{j}) \otimes$\\ 
$ F(X_{j+1} \otimes \ldots \otimes X_{j+p}) \otimes F(X_{j+p+1}) \otimes \ldots \otimes F(X_{n+p+i-1}) \otimes R $\\ $\otimes \ldots \otimes F(X_{m+n+p-2}) $ \end{tabular}} 
\arrow[d,"{\_ \otimes \ldots \otimes \rho^R(X_{i+n+p},\ldots, X_{m+n+p-2})}"] \\ \boxed{\begin{tabular}{l} $U \otimes F(X_1) \otimes...\otimes F(X_i) \otimes U \otimes F(X_{i+1}) \otimes... \otimes F(X_{j}) \otimes$\\ $\ldots \otimes F(X_{j+p}) \otimes F(X_{j+p+1}) \otimes \ldots \otimes F(X_{n+p+i-1}) \otimes R$\\ $\otimes \ldots F(X_{m+n+p-2}) $ \end{tabular} } 
\arrow[d,"{\_ \otimes \ldots \otimes g_{X_{i+1},\ldots,X_j,X_{j+1} \otimes \ldots \otimes X_{j+p}, \ldots, X_{n+p+i-1}} \otimes \_ \otimes \ldots}"] 
\& \& \& \& 
\boxed{\begin{tabular}{l} $U \otimes F(X_1) \otimes...\otimes F(X_i) \otimes$\\
$ F(X_{i+1} \otimes \ldots \otimes X_j \otimes \ldots \otimes X_{j+p} \otimes \ldots \otimes X_{i+n+p-1}) \otimes$\\ $\ldots \otimes F(X_{m+n+p-2}) \otimes R$ \end{tabular}} 
\arrow[d,"{f_{X_1,\ldots,X_i,X_{i+1}\otimes \ldots X_j \otimes \ldots \otimes X_{j+p} \otimes \ldots \otimes X_{i+n+p-1}, \ldots, X_{m+n+p-2}} \otimes \_ }"]\\ 
\boxed{ \begin{tabular}{l} $U \otimes F(X_1) \otimes...\otimes F(X_i) \otimes F(X_{i+1} \otimes \ldots \otimes X_{n+p+i-1}) \otimes R \otimes$\\ $R \otimes \ldots \otimes F(X_{m+n+p-2}) $ \end{tabular} } 
\arrow[d,"{\_ \otimes \ldots \otimes \_ \otimes \rho^R(X_{i+n+p},\ldots, X_{m+n+p-2})}"] 
\& \& \& \& 
\boxed{\begin{tabular}{l} $F(X_1) \otimes...\otimes F(X_i) \otimes$\\
$ F(X_{i+1} \otimes \ldots \otimes X_j \otimes \ldots \otimes X_{j+p} \otimes \ldots \otimes X_{i+n+p-1}) \otimes$\\ $\ldots \otimes F(X_{m+n+p-2}) \otimes R \otimes R$ \end{tabular}} 
\arrow[d, "{\ldots \otimes \mu }"] \\ 
\boxed{ \begin{tabular}{l} $U \otimes F(X_1) \otimes...\otimes F(X_i) \otimes F(X_{i+1} \otimes \ldots \otimes X_{n+p+i-1}) \otimes R \otimes$\\ $F(X_{i+n+p}) \otimes \ldots \otimes F(X_{m+n+p-2}) \otimes R$ \end{tabular} } 
\arrow[d,"{\_ \otimes \ldots \otimes \rho^R(X_{i+n+p},\ldots, X_{m+n+p-2}) \otimes \_ }"] 
\& \& \& \& 
\boxed{\begin{tabular}{l} $F(X_1) \otimes \ldots \otimes F(X_{m+n+p-2}) \otimes R$ \end{tabular}}\\
 \boxed{\begin{tabular}{l} $U \otimes F(X_1) \otimes...\otimes F(X_i) \otimes F(X_{i+1} \otimes \ldots \otimes X_{n+p+i-1})$\\ $F(X_{i+n+p}) \otimes \ldots \otimes F(X_{m+n+p-2}) \otimes R \otimes R$ \end{tabular} } 
 \arrow[d,"{\_ \otimes \ldots \otimes \mu }"] 
 \& \& \& \& \\
 \boxed{\begin{tabular}{l} $U \otimes F(X_1) \otimes \ldots \otimes F(X_i) \otimes F(X_{i+1} \otimes \ldots \otimes X_{n+p+i-1})$\\ $F(X_{i+n+p}) \otimes \ldots \otimes F(X_{m+n+p-2}) \otimes R$ \end{tabular} } 
 \arrow[rrrr, "{ f_{X_1,\ldots,X_i,X_{i+1}\otimes \ldots \otimes X_j\otimes \ldots \otimes X_{j+p} \otimes \ldots \otimes X_{i+n+p-1}, \ldots, X_{m+n+p-2}} \otimes \_ }"] 
 \& \& \& \& 
 \boxed{\begin{tabular}{l} $F(X_1) \otimes \ldots \otimes F(X_{m+n+p-2}) \otimes R \otimes R$ \arrow[uu,swap,"{ \ldots \otimes \mu}"] \end{tabular} } 
\end{tikzcd}
}
\end{center}

For the left hand side of \eqref{3.4tt}, we have
\begin{align}
&\big((f \diamond_i g) \diamond_{j} h)\big)_{X_1, \ldots ,X_{m+n+p-2}}= \nonumber\\ 
&(F(X_1) \otimes \ldots \otimes F(X_{m+n+p-2}) \otimes \mu)~ \circ 
(F(X_1) \otimes \ldots \otimes F(X_j) \otimes F(X_{j+1} \otimes \ldots \otimes X_{j+p}) \otimes F(X_{j+p+1}) \otimes \ldots \otimes F(X_{m+n+p-2}) \otimes \mu \otimes R)~\circ \nonumber\\
&(f_{X_1,...,X_i, X_{i+1} \otimes...X_j \otimes... \otimes X_{j+p}\otimes...\otimes X_{i+n+p-1}, \ldots, X_{m+n+p-2}} \otimes R)  \circ \nonumber\\
&(F(X_1) \otimes \ldots \otimes F(X_i) \otimes F(X_{i+1} \otimes \ldots \otimes X_j \otimes X_{j+1} \otimes \ldots \otimes X_{j+p} \otimes X_{j+p+1} \otimes \ldots \otimes X_{n+i+p-1}) \otimes \rho^R(X_{n+i+p}, \ldots, X_{m+n+p-2}) \otimes R) \circ \nonumber\\
&(U \otimes F(X_1) \otimes \ldots \otimes F(X_i) \otimes g_{X_{i+1},...,X_{j},X_{j+1}\otimes \ldots \otimes X_{j+p}, \ldots, X_{i+n+p-1}} \otimes F(X_{n+i+p}) \otimes \ldots \otimes F(X_{m+n+p-2}) \otimes R)~ \circ \label{3.5co}\\
&(U \otimes \rho^U(X_1,...,X_i) \otimes F(X_{i+1}) \otimes ...F(X_j) \otimes F(X_{j+1} \otimes \ldots \otimes X_{j+p}) \otimes F(X_{j+p+1}) \otimes \ldots \otimes F(X_{m+n+p-2}) \otimes R)~ \circ \nonumber\\
&(\Delta \otimes F(X_1) \otimes...\otimes F(X_j) \otimes F(X_{j+1} \otimes \ldots \otimes X_{j+p}) \otimes F(X_{j+p+1}) \otimes \ldots \otimes F(X_{m+n+p-2}) \otimes R)~ \circ \nonumber\\
&(U \otimes F(X_1) \otimes...\otimes F(X_j) \otimes F(X_{j+1} \otimes \ldots \otimes X_{j+p}) \otimes \rho^R(X_{j+p+1}, \ldots, X_{m+n+p-2})) \circ \nonumber\\
&(U \otimes F(X_1) \otimes...\otimes F(X_j) \otimes h_{X_{j+1} \otimes \ldots \otimes X_{j+p}} \otimes \ldots \otimes F(X_{m+n+p-2}))\circ \nonumber\\
&(U \otimes \rho^U(X_1,...,X_j) \otimes F(X_{j+1}) \otimes \ldots \otimes  F(X_{m+n+p-2})) \circ 
(\Delta \otimes F(X_1) \otimes...\otimes F(X_{m+n+p-2})) \nonumber
\end{align}
Using the associativity of $\mu$ and the coassociativity of $\Delta$, we see that the expression in \eqref{3.5co} is equal to the composition starting from the left vertical column of the previous commutative diagram. The result now follows.
\end{proof}

\begin{prop}\label{p3.4f}
Let $F:\mathscr C\longrightarrow \mathscr D$ be a   monoidal  functor between monoidal  categories. Let $(U,\rho^U,\Delta,\epsilon)$ be a coalgebra object  and $(R,\rho^R, \mu,\eta)$ be an algebra object in the centralizer $\mathcal Z(F)$ of $F$. Then,   $C^\bullet_{DY}(F,U, R)$ is a weak comp algebra.
\end{prop}
\begin{proof}
The condition {(1)} in Definition \ref{wg} is clear from the definition of the operations $\diamond_i$ on $C^\bullet_{DY}(F,U,R)$. The conditions {(2)} and {(3)} follow from Lemma \ref{con2}. The condition {(4)} may also be verified directly. 
\end{proof}

The next result shows that the operations $\cup$ and $\sqcup$ on the complex $C^\bullet_{DY}(F,U,R)$, as well as the differential $\partial$ can be recovered in terms of the operations $\diamond_i$ making $C^\bullet_{DY}(F,U,R)$ a weak comp algebra. 

\begin{lem}\label{cupsq}
Let $(U,\rho^U,\Delta,\epsilon)$ be a coalgebra object and $(R,\rho^R, \mu,\eta)$ be an algebra object in the category $\mathcal Z(F)$. Then, the following hold:
\begin{itemize}
\item[(a)] The operations $\cup$ and $\sqcup$ on $C^\bullet_{DY}(F,U,R)$ are also given by
\begin{equation*}
f \cup g=(\pi \diamond_0 f) \diamond_m g   \qquad \text{and} \qquad f \sqcup g=(\pi \diamond_1 g) \diamond_0 f
\end{equation*} for $f\in C^m_{DY}(F,U,R)$ and $g\in C^n_{DY}(F,U,R)$, where $m$, $n\geq 0$. 
\item[(b)] For any $m \geq 0$ and $f\in C^m_{DY}(F,U,R)$, the differential $\partial:C^m_{DY}(F,U,R) \longrightarrow C^{m+1}_{DY}(F,U,R)$ is also given by
$$\partial(f)=(-1)^{m-1} \pi \diamond_0 f-\sum\limits_{i=1}^m (-1)^{i-1} f \diamond_{i-1} \pi~ + \pi \diamond_1 f$$
\end{itemize}
\end{lem}
\begin{proof}
\emph{(a)} For any $X_1, \ldots, X_{m+n} \in Ob(\mathscr{C})$, we  have
\begin{align*}
&{((\pi \diamond_0 f) \diamond_m g)}_{X_1, \ldots, X_{m+n}}\\
&=  (F(X_1\otimes \ldots \otimes X_m) \otimes F(X_{m+1}\otimes \ldots \otimes X_{m+n}) \circ \mu) \circ 
(F(X_1\otimes \ldots \otimes X_m) \otimes F(X_{m+1}\otimes \ldots\otimes X_{m+n}) \otimes \mu \otimes R) \circ \\
& \quad F(X_1\otimes \ldots \otimes X_m) \otimes F(X_{m+1}\otimes \ldots \otimes X_{m+n}) \otimes \eta \otimes R \otimes R) \circ (\epsilon \otimes F(X_1\otimes \ldots \otimes X_m) \otimes F(X_{m+1}\otimes \ldots \otimes X_{m+n}) \otimes R \otimes R) \circ \\
& \quad (U \otimes F(X_1\otimes \ldots \otimes X_m) \otimes \rho^R(X_{m+1},\ldots, X_{m+n}) \otimes R ) \circ (U \otimes f_{X_1, \ldots, X_m} \otimes F(X_{m+1}\otimes \ldots \otimes X_{m+n}) \otimes R) \circ\\
&~~~~ (\Delta \otimes F(X_1) \otimes \ldots \otimes F(X_m) \otimes  F(X_{m+1}\otimes \ldots \otimes X_{m+n}) \otimes R) \circ (U \otimes F(X_1) \otimes \ldots \otimes F(X_m) \otimes  g_{X_{m+1}, \ldots, X_{m+n}}) \circ \\
&~~~~ (U\otimes \rho^U(X_1, \ldots, X_m)\otimes F(X_{m+1})\otimes \ldots \otimes F(X_{m+n})) \circ (\Delta\otimes F(X_1)\otimes \ldots \otimes F(X_{m+n}))\\
&=(F(X_1\otimes \ldots \otimes X_m) \otimes F(X_{m+1}\otimes \ldots \otimes X_{m+n}) \circ \mu) \circ (F(X_1)\otimes \ldots \otimes F(X_m) \otimes \rho^R(X_{m+1},\ldots,X_{m+n}) \otimes R) \circ\\ &~~~~(f_{X_1, \ldots, X_m} \otimes F(X_{m+1}\otimes \ldots \otimes X_{m+n}) \otimes R) \circ (\epsilon \otimes U \otimes F(X_1) \otimes \ldots \otimes F(X_m) \otimes F(X_{m+1}\otimes \ldots \otimes X_{m+n}) \otimes R) \circ\\
&~~~~ (\Delta \otimes F(X_1) \otimes \ldots \otimes F(X_m) \otimes  F(X_{m+1}\otimes \ldots \otimes X_{m+n} \otimes R)\circ (U \otimes F(X_1) \otimes \ldots \otimes F(X_m) \otimes  g_{X_{m+1}, \ldots, X_{m+n}}) \circ \\
\smallskip
&~~~~ (U\otimes \rho^U(X_1, \ldots, X_m)\otimes F(X_{m+1})\otimes \ldots \otimes F(X_{m+n})) \circ  (\Delta\otimes F(X_1)\otimes \ldots \otimes F(X_{m+n}))\\
&=(F(X_1\otimes \ldots \otimes X_m) \otimes F(X_{m+1}\otimes \ldots \otimes X_{m+n}) \circ \mu) \circ (F(X_1)\otimes \ldots \otimes F(X_m) \otimes \rho^R(X_{m+1},\ldots,X_{m+n}) ) \circ \\
&~~~~(f_{X_1, \ldots, X_m}\otimes g_{X_{m+1}, \ldots, X_{m+n}}) \circ (U\otimes \rho^U(X_1, \ldots, X_m)\otimes F(X_{m+1})\otimes \ldots\otimes F(X_{m+n}))  \circ\\
&~~~~(\Delta\otimes F(X_1)\otimes \ldots \otimes F(X_{m+n}))\\
&=(f\cup g)_{X_1, \ldots, X_{m+n}}
\end{align*}
where we use the unitality (resp. counitality) of $\eta$ (resp. $\epsilon$) in the second (resp. third) equality. The expression for $\sqcup$ may be verified similarly.


\medskip
(b) This may also be verified by direct computation.
\end{proof}

\begin{prop}\label{p3.6y}
Let $F:\mathscr C\longrightarrow \mathscr D$ be a  monoidal  functor between monoidal categories. Let $(U,\rho^U,\Delta,\epsilon)$ be a coalgebra object and $(R,\rho^R, \mu,\eta)$ be an algebra object in the centralizer $\mathcal Z(F)$ of $F$. Let $\eta \circ \epsilon \in C^0_{DY}(F,U,R)=Nat(U\otimes F^{\otimes 0},F^{\otimes 0} \otimes R)$ be the natural transformation given by
\begin{equation}
\eta \circ \epsilon : U\otimes F^{\otimes 0}(-)\xrightarrow{\qquad\qquad\epsilon\otimes F^{\otimes 0}(-)\qquad\qquad}  F^{\otimes 0}(-) \xrightarrow{\qquad \qquad F^{\otimes 0}(-) \otimes \eta \qquad \qquad} F^{\otimes 0} \otimes R 
\end{equation}
Then, we have the following
\begin{itemize}
\item[(a)] The complex $C^\bullet_{DY}(F,U,R)$ is a differential graded associative algebra with product $\cup$ and unit $\eta \circ \epsilon \in C^0_{DY}(F,U,R)$.
\item[(b)] The complex $C^\bullet_{DY}(F,U,R)$ is a differential graded associative algebra with product $\sqcup$ and unit $\eta \circ \epsilon \in C^0_{DY}(F,U,R)$.
\end{itemize}
\end{prop}
\begin{proof}
From Proposition \ref{p3.4f}, we know that the complex $C^\bullet_{DY}(F,U,R)$ is a weak comp algebra. From the description of 
$\cup$ and $\sqcup$ in Lemma \ref{cupsq}(a) and the properties of a weak comp algebra
obtained in  \cite[Proposition 4.5]{brz}, it now follows that $C^\bullet_{DY}(F,U,R)$ is a graded associative algebra for both the operations $\cup$ and 
$\sqcup$. 

\smallskip
Additionally, from the description of the differential $\partial$ in Lemma \ref{cupsq}(b) and the properties of a weak comp algebra
obtained in  \cite[Proposition 4.7]{brz},  it follows that both $(C^\bullet_{DY}(F,U,R),\cup, \partial^\bullet)$ and $(C^\bullet_{DY}(F,U,R),\sqcup, \partial^\bullet)$ are differential graded algebras. Moreover, for any $f \in C^m_{DY}(F,U,R)$ and $X_1, \ldots, X_m\in Ob(\mathscr C)$, using the counital and unital conditions in \eqref{2.8bv} and \eqref{Rcom}, we have
\begin{align*}
(f \cup (\eta \circ \epsilon))_{X_1, \ldots, X_m}&=(F(X_1) \otimes \ldots \otimes F(X_m) \otimes \mu) \circ (f_{X_1, \ldots, X_m} \otimes (\eta \circ\epsilon)) \circ (U \otimes \rho^U(X_1, \ldots, X_m)) \circ (\Delta \otimes F(X_1) \otimes \ldots \otimes F(X_m))\\
&=f_{X_1, \ldots, X_m} 
\end{align*}
Also, we can check that $((\eta \circ \epsilon) \cup f)_{X_1, \ldots, X_m}=f_{X_1, \ldots, X_m}$. This shows that $\eta \circ \epsilon$ is a unit for the product $\cup$ on $C^\bullet_{DY}(F,U,R)$. Similarly, it may  be verified that $\eta \circ \epsilon$ is also a unit for the product $\sqcup$.
\end{proof}

On its own, neither of the operations $\cup$ and $\sqcup$ on $C^\bullet_{DY}(F,U,R)$ is graded commutative. However, the following main result shows that we  have a version of ``graded commutativity'' involving the operations $\cup$ and $\sqcup$. 

\begin{Thm}\label{mdmsinhgr}
Let $F:\mathscr C\longrightarrow \mathscr D$ be a monoidal  functor between monoidal  categories. Let $(U,\rho^U,\Delta,\epsilon)$ be a coalgebra object and $(R,\rho^R, \mu,\eta)$ be an algebra object in the centralizer $\mathcal Z(F)$ of $F$. Then, 
\begin{itemize}
\item[(a)] There are cup products on the Davydov-Yetter cohomology $H^\bullet_{DY}(F,U,R)$:
\begin{align*}
\cup: H^m_{DY}(F,U,R) \otimes H^n_{DY}(F,U,R) \longrightarrow H^{m+n}_{DY}(F,U,R)\qquad\qquad 
\sqcup: H^m_{DY}(F,U,R) \otimes H^n_{DY}(F,U,R) \longrightarrow H^{m+n}_{DY}(F,U,R)\quad m,n\geq 0
\end{align*}
making $H^\bullet_{DY}(F,U,R)$ into a graded associative algebra.

\smallskip
\item[(b)] For cohomology classes $\bar{f} \in H^m_{DY}(F,U,R)$ and $\bar{g} \in H^n_{DY}(F,U,R)$,  we have
\begin{equation}\label{bnjr}
\bar{f} \cup \bar{g} = (-1)^{mn}~ \bar{g} \sqcup \bar{f}
\end{equation}
\end{itemize}
\end{Thm}
\begin{proof}
From Proposition \ref{p3.6y}, we know that there are cup products $\cup$ and $\sqcup$ on the complex $C_{DY}^\bullet(F,U,R)$, making it an associative algebra. From Lemma \ref{L2.1} and 
Lemma \ref{L2.2}, it follows that the products $\cup$ and $\sqcup$ descend to the level of cohomology. This proves (a). Since $C^\bullet_{DY}(F,U,R)$
is a weak comp algebra, the result of (b) now follows by applying \cite[Corollary 4.9]{brz}.
\end{proof}

\begin{eg}\label{grat}
\emph{Let  $\mathscr C=\mathscr D$ and $F:=id_{\mathscr C}:\mathscr C\longrightarrow \mathscr C$. In that case, the centralizer $\mathcal Z(F)$ reduces to the Drinfeld center of the monoidal  category $\mathscr C$, i.e., its objects are pairs $(U,\rho^U)$, where (see \cite[$\S$ 7.13]{Etbook})
\begin{equation*}
\rho^U=\{\rho^U(X):U\otimes X\xrightarrow{\quad\cong\quad } X\otimes U\}_{X\in Ob(\mathscr C)}
\end{equation*} consists of natural isomorphisms which satisfy compatibility conditions as in \eqref{2.2rc}. As mentioned in Section 2, the centralizer 
$\mathcal Z(F)$ always carries a monoidal structure. Accordingly, if $(U,\rho^U,\Delta,\epsilon)$ is a coalgebra object and $(R,\rho^R,\mu,\eta)$ is an algebra object in the Drinfeld center
$\mathcal Z(id_{\mathscr C})$ of $\mathscr C$, it follows that $C^\bullet_{DY}(\mathscr C,U,R):=C^\bullet_{DY}(id_{\mathscr C},U,R)$ is a weak comp algebra and the Davydov-Yetter cohomology $H^\bullet_{DY}(\mathscr C,U,R):=H^\bullet_{DY}(id_{\mathscr C},U,R)$ is equipped with cup products $\cup$ and $\sqcup$ which are related to each other as in \eqref{bnjr}. We also recall that when $U=R=1$ is the unit object in $\mathscr C$, the groups $H^\bullet_{DY}(\mathscr C):=H^\bullet_{DY}(id_{\mathscr C},1,1)$ are known as the Davydov-Yetter cohomology groups of
$\mathscr C$ (see \cite[$\S$ 7.22.1]{Etbook}).}
\end{eg}

\begin{eg}\label{eg3.8}
\emph{We recall that a monoidal category $\mathscr D$ is said to be a tensor category if it is also rigid, i.e.,  every object $V\in \mathscr D$ admits a left dual $^\vee V$ and a right dual 
$V^\vee$ (see, for instance, \cite[$\S$ 2.10]{Etbook}).  Further, a $k$-linear tensor category is said to be finite if it is equivalent to the category of finite dimensional representations of a finite dimensional $k$-algebra (see, for instance, \cite[$\S$ 3.1]{Gai}).  Let $F:\mathscr C\longrightarrow \mathscr D$ be a tensor functor between finite tensor categories $\mathscr C$ and $\mathscr D$. Additionally, suppose that $F$ is such that for each $V\in \mathscr D$, the object  
\begin{equation}\label{centmon}
Z_F(V):=\int^{X\in Ob(\mathscr C)} F(X)^\vee\otimes V\otimes F(X)
\end{equation}
exists. This always happens, for instance, when $F$ is an exact functor (see \cite[$\S$ 3.3]{Gai}). In that case, the association $V\mapsto Z_F(V)$ defines a monad on $\mathscr D$, which is known as the central monad of $F$. Then, $Z_F$ is a Hopf monad in the sense of \cite{BV}, and the centralizer 
$\mathcal Z(F)$ of $F$ is isomorphic as a tensor category to the Eilenberg-Moore category $Z_F-mod$ of modules over the monad $Z_F$ (see \cite[Proposition 3.10]{Gai}). Accordingly, if $(U,\rho^U,\Delta,\epsilon)$ is a coalgebra object and $(R,\rho^R,\mu,\eta)$ is an algebra object in 
$Z_F-mod$, it follows that the Davydov-Yetter complex $C^\bullet_{DY}(F,U,R)$ with coefficients in $U$ and $R$ is a weak comp algebra and the cohomology $H^\bullet_{DY}(F,U,R)$ is equipped with cup products $\cup$ and $\sqcup$ which are related to each other as in \eqref{bnjr}.  For instance, let $H$ be a finite dimensional Hopf algebra over $k$, and let $\mathscr C:=H-mod$ be the category of left modules over $H$ that are also finite dimensional over $k$. Then, $H-mod$ is a finite tensor category, and we can take $F$ to be the forgetful functor from $H-mod$ to $Vect_k$ (see \cite[$\S$ 1]{Gai}).}
\end{eg}

\begin{eg}
\emph{Let $G$ be a group. We consider Example \ref{grat} with $\mathscr C=Vec_G$, the monoidal  category of $G$-graded finite dimensional $k$-vector spaces. Accordingly, every object 
$V\in \mathscr C=Vec_G$ is of the form $V={\oplus_{g \in G}} V_g$. By \cite[$\S$ 8.5.4]{Etbook}, the Drinfeld center of $Vec_G$ consists of $G$-equivariant objects of $Vec_G$ with respect to the conjugation action of $G$ on $Vec_G$, i.e., those $V\in Vec_G$ for which $V_{gxg^{-1}}\cong V_x$ for all $x$, $g\in G$. Then, the Drinfeld center $\mathcal Z(Vec_G)$ has a monoidal structure inherited from $Vec_G$. If $(U,\rho^U,\Delta,\epsilon)$ is a coalgebra object and $(R,\rho^R,\mu,\eta)$ is an algebra object in 
$\mathcal Z(Vec_G)$, it follows that the Davydov-Yetter complex $C^\bullet_{DY}(Vec_G,U,R)$ with coefficients in $U$ and $R$ is a weak comp algebra and the cohomology $H^\bullet_{DY}(Vec_G,U,R)$ is equipped with cup products $\cup$ and $\sqcup$ which are related to each other as in \eqref{bnjr}.}
\end{eg}

\begin{eg}
\emph{Let $H$ be a finite dimensional Hopf algebra over a field $k$, with invertible antipode $S$. We consider Example \ref{grat}  with $\mathscr C=H-mod$, the monoidal  category of left modules over $H$ that are also finite dimensional over $k$. In that case, we know  that the Drinfeld center $\mathcal Z(H-mod)$ can be identified as a monoidal  category with $D(H)-mod$, the category of finite dimensional representations of the quantum double $D(H)$ of $H$ (see \cite[$\S$ 8.5.6]{Etbook}). Let $(U,\rho^U,\Delta,\epsilon)$ be a coalgebra object and $(R,\rho^R,\mu,\eta)$ be an algebra object in the category $D(H)-mod$. It follows that the Davydov-Yetter complex $C^\bullet_{DY}(H-mod,U,R)$ with coefficients in $U$ and $R$  is a weak comp algebra and the cohomology $H^\bullet_{DY}(H-mod,U,R)$ is equipped with cup products $\cup$ and $\sqcup$ which are related to each other as in \eqref{bnjr}. More explicitly,  we know from   \cite[Lemma 5.2]{Fai} that the Davydov-Yetter complex 
\begin{equation} \label{3.9mbl} 
\begin{array}{c} C^\bullet_{DY}(H-mod,U,R) \cong  \text{Hom}_H\big(U, (H^{\otimes \bullet} \otimes R)_{\text{ad}}\big)\qquad\qquad 
f\mapsto (u\mapsto f_{H, \ldots, H}(u\otimes 1_H\otimes \ldots \otimes 1_H))
\\
\end{array} \end{equation}
where we take  $(H^{\otimes 0} \otimes R)_{\text{ad}}=R$ and for \(n \geq 1\), the space
$\left(H^{\otimes n} \otimes R\right)_{\mathrm{ad}}$ is the left $H$-module
with the action
\begin{equation}
\begin{aligned}
h \cdot (x_1 \otimes \cdots \otimes x_n \otimes r)
={}& h_{(1)}x_1S\bigl(h_{(2n+1)}\bigr)
\otimes \cdots \otimes
h_{(n)}x_nS\bigl(h_{(n+2)}\bigr) \otimes h_{(n+1)} \cdot r
\end{aligned}
\end{equation}
for \(x_i \in H\),  \(1 \leq i \leq n\),  \(h \in H\),  and \(r \in R\). The differential on the complex $\text{Hom}_H\big(U, (H^{\otimes \bullet} \otimes R)_{\text{ad}}\big)$  is given by
\begin{equation}
\begin{aligned}
\delta^n(\varphi)
={}&
(H \otimes \varphi) \circ \lambda_U
+ \sum_{i=1}^{n} (-1)^i
\left(
H^{\otimes(i-1)}
\otimes \Delta_H
\otimes H^{\otimes(n-i)}
\otimes R
\right) \circ \varphi 
+ (-1)^{n+1}
\left(
H^{\otimes n} \otimes \lambda_R
\right) \circ \varphi.
\end{aligned}
\end{equation}
for $\varphi \in  \mathrm{Hom}_H\big(U, (H^{\otimes n} \otimes R)_{\mathrm{ad}}\big)$, where $\lambda_U$ and $\lambda_R$ are respectively the left $H$-coactions on $U, R \in D(H)-{mod} \simeq {}_H^H\mathcal{YD}$ (the category of left-left Yetter-Drinfeld modules over $H$) (see \cite[Theorem IX.5.2]{ka}) and $\Delta_H$ denotes the comultiplication map on $H$.}

\smallskip
\emph{For any $u \in U, ~m \geq 0$, $f \in C^m_{DY}(H-mod,U,R) $ and $x_i \in X_i$, $X_i\in H-mod$, $1 \leq i \leq m$,  we write $f_{X_1, \ldots, X_m}(u \otimes x_1 \otimes \ldots \otimes x_m)=f^{x_1}_u \otimes \ldots \otimes f^{x_m}_u \otimes f^R_u\in X_1\otimes ...\otimes X_m\otimes R$, with summation omitted. We also write $\rho^U(H)(u \otimes x)=x_\alpha \otimes u^\alpha \in H \otimes U$ and $\rho^R(H)(r \otimes x)=x_\beta \otimes r^\beta \in H \otimes R$, with summation omitted.
Using the isomorphism in \eqref{3.9mbl}, the cup products $\cup$ and $\sqcup$ on the complex $C^\bullet_{DY}(H-mod,U,R) $ may now be expressed explicitly as
\begin{equation*}
(f \cup g)(u)=f^{{1_H}_{\alpha_1}}_{u_{(1)}} \otimes \ldots \otimes f^{{1_H}_{\alpha_m}}_{u_{(1)}}  \otimes \big(g^{1_H}_{u_{(2)}^{\alpha_1 \cdots \alpha_m}}\big)_{\beta_1}  \otimes \ldots \otimes \big(g^{1_H}_{u_{(2)}^{\alpha_1 \cdots \alpha_m}}\big)_{\beta_n}  \otimes \mu((f^R_{u_{(1)}})^{\beta_1\cdots \beta_n}  \otimes g^R_{u_{(2)}^{\alpha_1 \cdots \alpha_m}} )
\end{equation*}
\begin{equation*}
(f \sqcup g)(u)=(f^{1_H}_{u_{(2)}})_{\alpha_1} \otimes \ldots \otimes (f^{1_H}_{u_{(2)}})_{\alpha_m} \otimes g^{{1_H}_{\beta_1}}_{u_{(1)}^{\alpha_1...\alpha_m}} \otimes \ldots g^{{1_H}_{\beta_n}}_{u_{(1)}^{\alpha_1...\alpha_m}} \otimes \mu\big(g^R_{u_{(1)}^{\alpha_1...\alpha_m}} \otimes (f^R_{u_{(2)}})^{\beta_1...\beta_n}\big) 
\end{equation*}
for $f \in C^m_{DY}(H-mod,U,R) \cong \text{Hom}_H\big(U, (H^{\otimes m} \otimes R)_{\text{ad}}\big),$ $g \in C^n_{DY}(H-mod,U,R) \cong  \text{Hom}_H\big(U, (H^{\otimes n} \otimes R)_{\text{ad}}\big)$ and $u \in U$, with $\Delta(u)=u_{(1)} \otimes u_{(2)}$ in Sweedler notation with summation sign omitted. In the same way, the  operations $(f \diamond_i g)$ for $0 \leq i <m$ may also be explicitly described as follows
\begin{equation*}
(f \diamond_i g)(u)=f^{{1_H}_{\alpha_1}}_{u_{(1)}} \otimes \ldots \otimes f^{{1_H}_{\alpha_i}}_{u_{(1)}}  \otimes \ldots \otimes f_{u_{(1)}}^{g^{1_H}_{{u_{(2)}}^{\alpha_1 \cdots \alpha_i}}\otimes ...\otimes g^{1_H}_{{u_{(2)}}^{\alpha_1 \cdots \alpha_i}}} \otimes f^{{1_H}_{\beta_1}}_{u_{(1)}} \otimes \ldots \otimes f^{{1_H}_{\beta_{m-i-1}}}_{u_{(1)}} \otimes \ldots \otimes \mu\big(f^R_{u_{(1)}} \otimes (g^R_{{u_{(2)}}^{\alpha_1 \cdots \alpha_i}})^{\beta_1...\beta_{m-i-1}}\big)
\end{equation*}}
\end{eg}

\section{Gerstenhaber algebra structure}
We continue with $(U,\rho^U,\Delta,\epsilon)$ being a coalgebra object and $(R, \rho^R,\mu,\eta)$ being an algebra object in the centralizer $\mathcal Z(F)$ of the monoidal functor $F:
\mathscr C\longrightarrow \mathscr D$.  In Theorem \ref{mdmsinhgr}, we showed that the cup products $\cup$ and $\sqcup$ on $H^\bullet_{DY}(F,U,R)$ are related by a version of 
``graded commutativity'' given by
\begin{equation}\label{bnjr4c}
\bar{f} \cup \bar{g} = (-1)^{mn}~ \bar{g} \sqcup \bar{f}
\end{equation} for cohomology classes $\bar{f} \in H^m_{DY}(F,U,R)$ and $\bar{g} \in H^n_{DY}(F,U,R)$, $m$, $n\geq 0$. 
Our final aim in this paper is to restrict to subcomplexes of $C^\bullet_{DY}(F,U,R)$ such that the cup product $\cup$ on cohomology is graded commutative, and whose cohomology is a Gerstenhaber algebra.  

\smallskip Building on the analogy between half-braidings and entwining structures, we will now construct a counterpart of the equivariant complex of Brzezi{\'n}ski \cite[$\S$ 5]{brz} in the Davydov-Yetter context. For this, we start with a morphism $\,\overline{\!\Psi\!}\,: U \otimes R \longrightarrow R \otimes U$  in the category  $\mathcal Z(F)$. This means  that the following diagram commutes:

\begin{equation}\label{morp}
\begin{tikzcd}[column sep=large, row sep=large]
U \otimes R \otimes F(X)
  \arrow[r, "U \otimes \rho^R(X)"]
  \arrow[d, "\,\overline{\!\Psi\!}\,\otimes {F(X)}"]
&
U \otimes F(X) \otimes R
  \arrow[r, "\rho^U(X) \otimes R"]
&
F(X) \otimes U \otimes R
  \arrow[d, "{F(X)} \otimes \overline{\Psi}"]
\\
R \otimes U \otimes F(X)
  \arrow[r, "R \otimes \rho^U(X)"]
&
R \otimes F(X) \otimes U
  \arrow[r, "\rho^R(X) \otimes U"]
&
F(X) \otimes R \otimes U
\end{tikzcd}
\end{equation}

We also suppose that the morphism $\,\overline{\!\Psi\!}\,$ is well behaved with respect to the multiplication $\mu$ and the unit $\eta$ of the algebra object $R$, i.e. the following diagrams commute
\begin{equation}\label{mult}
\begin{tikzcd}
U \otimes R \otimes R
\arrow[r, "{\,\overline{\!\Psi\!}\,\otimes R}"]
\arrow[d, swap, "{U \otimes \mu}"]
&
R \otimes U \otimes R
\arrow[r, "{R \otimes \,\overline{\!\Psi\!}\,}"]
&
R \otimes R \otimes U
\arrow[d, "{\mu \otimes U}"]\\
U \otimes R
\arrow[rr,"{\,\overline{\!\Psi\!}\,}"]
& 
&
R \otimes U 
\end{tikzcd} \qquad \qquad
\begin{tikzcd}
& U 
\arrow[ld, swap, "{U \otimes \eta }"]
\arrow[rd, "{\eta \otimes U}"]
&\\
U \otimes R 
\arrow[rr, "{\,\overline{\!\Psi\!}\,}"] 
& & R \otimes U
\end{tikzcd}
\end{equation}


%
%

For each $m \geq 0$, we  set $P^m(F,U,R,\,\overline{\!\Psi\!}\,)\subseteq C^m_{DY}(F,U,R)=Nat(U\otimes F^{\otimes m},F^{\otimes m}\otimes R)$ to be the collection of those natural tranformations $f:U \otimes F^{\otimes m} \longrightarrow F^{\otimes m} \otimes R$ for which the following diagram commutes
\begin{equation}\label{3.10cx}\small 
\begin{tikzcd}[column sep=2.4cm, row sep=large]
U \otimes F(X_1) \otimes \ldots \otimes F(X_m)
\arrow[r, "{\Delta \otimes F(X_1) \otimes \ldots \otimes F(X_m)}"]
&
U \otimes U \otimes F(X_1) \otimes \ldots \otimes F(X_m)
\arrow[r, "{U \otimes \rho^U(X_1,\ldots,X_m)}"]
\arrow[d, swap, "{ U \otimes f_{X_1,\ldots,X_m}}"]
&
U \otimes F(X_1) \otimes \ldots \otimes F(X_m) \otimes U
\arrow[dd,"{f_{X_1,\ldots,X_m} \otimes U }"]
\\
& U \otimes F(X_1) \otimes \ldots \otimes F(X_m) \otimes R
\arrow[d, swap, "{\rho^U(X_1,\ldots,X_m) \otimes R}"]
\\
&
F(X_1) \otimes \ldots \otimes F(X_m) \otimes U \otimes R
\arrow[r, "{F(X_1) \otimes \ldots \otimes F(X_m) \otimes \,\overline{\!\Psi\!}\,}"]
&
F(X_1) \otimes \ldots \otimes F(X_m) \otimes R \otimes U
 \end{tikzcd}
 \end{equation} for all $X_1$,...,$X_m\in Ob(\mathscr C)$. We now have the following result. 

\begin{prop}\label{ent11}
The subspaces $P^m(F,U,R,\,\overline{\!\Psi\!}\,)\subseteq C^m_{DY}(F,U,R)$ for $m\geq 0$ determine a subcomplex of $(C^\bullet_{DY}(F,U,R),\partial^\bullet)$.
\end{prop}
\begin{proof} We need to verify that $\partial(P^m(F,U,R,\,\overline{\!\Psi\!}\,))\subseteq P^{m+1}(F,U,R,\,\overline{\!\Psi\!}\,)$ for each $m\geq 0$. 
Let $m >0$ and $f\in P^m(F,U,R,\,\overline{\!\Psi\!}\,)$.  For $1 \leq i \leq m$ and $X_0,...,X_{m} \in \mathscr C$, it may be verified that direct computation that
\begin{align}\label{4.5mbsk}
&\big((\partial_if)_{X_0, \ldots, X_{m}}  \otimes U\big) \circ \big(U \otimes \rho^U(X_0, \ldots, X_{m})\big) \circ \big(\Delta \otimes F(X_0) \otimes \ldots \otimes F(X_m)  \big)=\notag \\
&\qquad \qquad \qquad \big(F(X_0) \otimes \ldots \otimes F(X_m) \otimes \,\overline{\!\Psi\!}\,\big) \circ \big(\rho^U(X_0, \ldots, X_{m}) \otimes R\big) \circ \big(U \otimes (\partial_if)_{X_0, \ldots, X_{m}} \big) \circ \big(\Delta \otimes F(X_0) \otimes \ldots \otimes F(X_m)  \big)
\end{align}
Further, we see that 
\begin{align} \label{4.6mbsk}
&\big(F(X_0) \otimes \ldots \otimes F(X_m) \otimes \,\overline{\!\Psi\!}\,\big) \circ \big(\rho^U(X_0, \ldots, X_{m}) \otimes R\big) \circ \big(U \otimes (\partial_{m+1}f)_{X_0,\ldots, X_{m}} \big) \circ \big(\Delta \otimes F(X_0) \otimes \ldots \otimes F(X_m)  \big)\notag \\
&~ = \big(F(X_0) \otimes \ldots \otimes F(X_m) \otimes \,\overline{\!\Psi\!}\,\big) \circ \big(\rho^U(X_0,\ldots, X_{m}) \otimes R\big) \circ \big(U \otimes F(X_0) \otimes \ldots \otimes F(X_{m-1}) \otimes \rho^R(X_m) \big) \circ \big(U \otimes f_{X_0,\ldots,X_{m-1}} \otimes F(X_m)\big) \circ\notag \\
& \qquad \big(\Delta \otimes F(X_0) \otimes \ldots \otimes F(X_m)  \big)\notag \\
& ~=  \big(F(X_0) \otimes \ldots \otimes F(X_{m-1}) \otimes \rho^R(X_m) \otimes U\big) \circ \big(F(X_0) \otimes \ldots \otimes F(X_{m-1}) \otimes R \otimes \rho^U(X_m)\big) \circ \big(F(X_0) \otimes \ldots \otimes F(X_{m-1}) \otimes \,\overline{\!\Psi\!}\,\otimes F(X_m) \big) \circ\notag \\
&  \qquad \big(\rho^U(X_0, \ldots, X_{m-1}) \otimes R \otimes F(X_m)\big) \circ \big(U \otimes f_{X_0,\ldots,X_{m-1}} \otimes F(X_m)\big) \circ \big(\Delta \otimes F(X_0) \otimes \ldots \otimes F(X_m)  \big)\notag \\
& \stackrel{\eqref{3.10cx}}{=}\big(F(X_0) \otimes \ldots \otimes F(X_{m-1}) \otimes \rho^R(X_m) \otimes U\big) \circ \big(F(X_0) \otimes \ldots \otimes F(X_{m-1}) \otimes R \otimes \rho^U(X_m)\big) \circ \big(f_{X_0,\ldots,X_{m-1}} \otimes U \otimes F(X_m)\big) \circ \notag \\
&  \qquad \big(U \otimes \rho^U(X_0, \ldots, X_{m-1}) \otimes F(X_m)\big) \circ \big(\Delta \otimes F(X_0) \otimes \ldots \otimes F(X_m)  \big)\notag \\
&~ =\big(F(X_0) \otimes \ldots \otimes F(X_{m-1}) \otimes \rho^R(X_m) \otimes U\big) \circ \big(f_{X_0,\ldots,X_{m-1}} \otimes F(X_m) \otimes U \big) \circ \big(U \otimes F(X_0) \otimes \ldots \otimes F(X_{m-1}) \otimes \rho^U(X_m)\big) \circ \notag \\
& \qquad \big(U \otimes \rho^U(X_0, \ldots, X_{m-1}) \otimes F(X_m)\big) \circ  \big(\Delta \otimes F(X_0) \otimes \ldots \otimes F(X_m)  \big)\notag \\
&~= \big((\partial_{m+1}f)_{X_0,\ldots, X_{m}} \otimes U\big) \big(U \otimes \rho^U(X_0, \ldots, X_{m}) \big) \circ  \big(\Delta \otimes F(X_0) \otimes \ldots \otimes F(X_m)  \big)
\end{align}
Also, we have
\begin{align}\label{4.7mbsk}
&\big(F(X_0) \otimes \ldots \otimes F(X_m) \otimes \,\overline{\!\Psi\!}\,\big) \circ \big(\rho^U(X_0, \ldots, X_{m}) \otimes R\big) \circ \big(U \otimes (\partial_{0}f)_{X_0,\ldots, X_{m}} \big) \circ \big(\Delta \otimes F(X_0) \otimes \ldots \otimes F(X_m)  \big)\notag \\
&~ =\big(F(X_0) \otimes \ldots \otimes F(X_m) \otimes \,\overline{\!\Psi\!}\,\big) \circ \big(\rho^U(X_0, \ldots, X_{m}) \otimes R\big) \circ  \big(U \otimes F(X_0) \otimes f_{X_1,\ldots,X_{m}}\big) \circ \notag \\
& \qquad  \big(U \otimes \rho^U(X_0) \otimes F(X_1) \otimes \ldots \otimes F(X_m)\big) \circ \big(\Delta \otimes F(X_0) \otimes \ldots \otimes F(X_m)  \big)\notag \\
& ~=\big(F(X_0) \otimes \ldots \otimes F(X_m) \otimes \,\overline{\!\Psi\!}\,\big) \circ \big(F(X_0) \otimes \rho^U(X_1,\ldots,X_{m}) \otimes R\big) \circ \big(F(X_0) \otimes U \otimes f_{X_1,\ldots,X_{m}} \big)\circ\notag \\
& \qquad \big(\rho^U(X_0) \otimes U \otimes F(X_1) \otimes \ldots \otimes F(X_m)\big) \circ \big(U \otimes \rho^U(X_0)  \otimes F(X_1) \otimes \ldots \otimes F(X_m)\big) \circ \big(\Delta \otimes F(X_0) \otimes \ldots \otimes F(X_m)  \big)\notag \\
&\stackrel{\eqref{2.8bv}}{=}\big(F(X_0) \otimes \ldots \otimes F(X_m) \otimes \,\overline{\!\Psi\!}\,\big) \circ \big(F(X_0) \otimes \rho^U(X_1,\ldots,X_{m}) \otimes R\big) \circ \big(F(X_0) \otimes U \otimes f_{X_1,\ldots,X_{m}} \big)\circ\notag \\
& \qquad \big(F(X_0) \otimes \Delta \otimes F(X_1) \otimes \ldots \otimes F(X_m) \big) \circ \big(\rho^U(X_0) \otimes F(X_1) \otimes \ldots \otimes F(X_m) \big)\notag \\
&\stackrel{\eqref{3.10cx}}{=}\big(F(X_0) \otimes f_{X_1,\ldots,X_{m}} \otimes U\big) \circ \big(F(X_0) \otimes U \otimes \rho^U(X_1, \ldots, X_{m})\big) \circ \big(F(X_0) \otimes \Delta \otimes F(X_1) \otimes \ldots \otimes F(X_m) \big) \circ\notag \\
& \qquad \big(\rho^U(X_0) \otimes F(X_1) \otimes \ldots \otimes F(X_m) \big)\notag \\
&\stackrel{\eqref{2.8bv}}{=}\big(F(X_0) \otimes f_{X_1,\ldots,X_{m}} \otimes U\big) \circ \big(F(X_0) \otimes U \otimes \rho^U(X_1, \ldots, X_{m})\big) \circ \big(\rho^U(X_0) \otimes U \otimes F(X_1) \otimes \ldots \otimes F(X_m)\big) \circ\notag \\
& \qquad \big(U \otimes \rho^U(X_0)  \otimes F(X_1) \otimes \ldots \otimes F(X_m)\big) \circ \big(\Delta \otimes F(X_0) \otimes \ldots \otimes F(X_m)  \big)\notag \\
&~= \big((\partial_{0}f)_{X_0,\ldots, X_{m}} \otimes U\big) \circ \big(U \otimes F(X_0) \otimes \rho^U(X_1, \ldots, X_{m})\big) \circ \big(U \otimes \rho^U(X_0)  \otimes F(X_1) \otimes \ldots \otimes F(X_m)\big) \circ\notag \\
&\qquad \big(\Delta \otimes F(X_0) \otimes \ldots \otimes F(X_m)  \big)\notag \\
&~=\big((\partial_{0}f)_{X_0,\ldots, X_{m}} \otimes U\big) \circ \big(U \otimes \rho^U(X_0, \ldots, X_{m})\big) \circ \big(\Delta \otimes F(X_0) \otimes \ldots \otimes F(X_m)  \big)
\end{align} From \eqref{4.5mbsk}, \eqref{4.6mbsk} and \eqref{4.7mbsk}, it follows that  $\partial_k(P^m(F,U,R,\,\overline{\!\Psi\!}\,))\subseteq P^{m+1}(F,U,R,\,\overline{\!\Psi\!}\,)$
for each $0\leq k\leq m+1$. This shows that $\partial(P^m(F,U,R,\,\overline{\!\Psi\!}\,))\subseteq P^{m+1}(F,U,R,\,\overline{\!\Psi\!}\,)$ for $m>0$.  
For $m=0$, it may be verified similarly that
\begin{align*}
(F(X_0) \otimes \,\overline{\!\Psi\!}\,) \circ (\rho^U(X_0) \otimes R) \circ (U \otimes (\partial^0f)_{X_0}) \circ (\Delta \otimes F(X_0))=((\partial^0f)_{X_0}) \otimes U) \circ (U \otimes \rho^U(X_0))\circ (\Delta \otimes F(X_0))
\end{align*}
This proves the result.
\end{proof}

Next, we recall the operations $\diamond_i$ on $C^\bullet_{DY}(F,U,R)$ as defined in \eqref{3.1q} and show that they  restrict to the subcomplex $\big(P^\bullet(F,U,R,\,\overline{\!\Psi\!}\,),\partial^\bullet\big)$. For any $f \in C^m_{DY}(F,U,R)$ and $g \in  C^n_{DY}(F,U,R)$ and $0 \leq i < m$, we have
$f\diamond_ig\in C^{m+n-1}_{DY}(F,U,R)$ defined by 
\begin{align}
(f \diamond_i g)_{X_1, \ldots, X_{m+n-1}}&=\big( F(X_1) \otimes \ldots \otimes F(X_{m+n-1}) \otimes \mu \big) \circ (f_{X_1, \ldots, X_i, X_{i+1} \otimes \ldots \otimes X_{i+n}, \ldots, X_{m+n-1}} \otimes R) \circ\notag\\
& \quad \big(U \otimes  F(X_1) \otimes \ldots \otimes F(X_i) \otimes F(X_{i+1} \otimes \ldots \otimes X_{i+n}) \otimes \rho^R(X_{i+n+1, \ldots, X_{m+n-1}}) \big) \circ\notag \\
&  \quad  (U \otimes F(X_1) \otimes \ldots \otimes F(X_i) \otimes g_{X_{i+1}, \ldots, X_{i+n}} \otimes F(X_{i+n+1})\otimes \ldots \otimes  F(X_{m+n-1}))  \circ \notag\\
& \quad \big((U \otimes \rho^U(X_1, \ldots, X_i)) \otimes F(X_{i+1}) \otimes \ldots \otimes  F(X_{m+n-1})\big) \circ  (\Delta \otimes F(X_1) \otimes \ldots\otimes F(X_{m+n-1}))
\end{align}
for $X_1, \ldots, X_{m+n-1}\in Ob(\mathscr C)$. For all other values of $i$, we have $f \diamond_i g=0$.  

\begin{lem}\label{closed}
For $f \in P^m(F,U,R,\,\overline{\!\Psi\!}\,)$ and $g \in P^n(F,U,R,\,\overline{\!\Psi\!}\,)$, we have each $f \diamond_i g \in P^{m+n-1}(F,U,R,\,\overline{\!\Psi\!}\,)$.
\end{lem}
\begin{proof}
Let $X_1, \ldots, X_{m+n-1}\in Ob(\mathscr C)$. For $0 \leq i < m$, we observe that
\begin{align*}
&\big((f \diamond_i g)_{X_1, \ldots, X_{m+n-1}} \otimes U\big) \circ \big(U \otimes \rho^U(X_1, \ldots, X_{m+n-1})\big) \circ  \big(\Delta \otimes F(X_1) \otimes \ldots \otimes F(X_{m+n-1})\big)\\
&=\big( F(X_1) \otimes \ldots \otimes F(X_{m+n-1}) \otimes \mu \otimes U\big) \circ \big(f_{X_1, \ldots, X_i, X_{i+1} \otimes \ldots \otimes X_{i+n}, \ldots, X_{m+n-1}} \otimes R \otimes U\big) \circ\\
&\quad \big(U \otimes  F(X_1) \otimes \ldots \otimes F(X_i) \otimes F(X_{i+1} \otimes \ldots \otimes X_{i+n}) \otimes \rho^R(X_{i+n+1}, \ldots, X_{m+n-1}) \otimes U\big) \circ \\
& \quad  \big(U \otimes F(X_1) \otimes \ldots \otimes F(X_i) \otimes g_{X_{i+1}, \ldots, X_{i+n}} \otimes F(X_{i+n+1})\otimes \ldots \otimes  F(X_{m+n-1}) \otimes U \big)  \circ \\
&  \quad \big(U \otimes \rho^U(X_1, \ldots, X_i) \otimes F(X_{i+1}) \otimes \ldots \otimes  F(X_{m+n-1}) \otimes U\big) \circ  (\Delta \otimes F(X_1) \otimes \ldots\otimes F(X_{m+n-1}) \otimes U) \circ \\
& \quad \big(U \otimes \rho^U(X_1, \ldots, X_{m+n-1})\big) \circ  \big(\Delta \otimes F(X_1) \otimes \ldots \otimes F(X_{m+n-1})\big)\\
&\stackrel{\eqref{2.8bv}}{=}\big( F(X_1) \otimes \ldots \otimes F(X_{m+n-1}) \otimes \mu \otimes U\big) \circ \big(f_{X_1, \ldots, X_i, X_{i+1} \otimes \ldots \otimes X_{i+n}, \ldots, X_{m+n-1}} \otimes R \otimes U\big) \circ
\\
& \quad \big(U \otimes F(X_1) \otimes \ldots \otimes F(X_i) \otimes F(X_{i+1} \otimes \ldots \otimes X_{i+n})   \otimes \rho^R(X_{i+n+1}, \ldots, X_{m+n-1})\otimes U\big) \circ\\
& \quad \big(U \otimes F(X_1) \otimes \ldots \otimes F(X_i) \otimes F(X_{i+1} \otimes \ldots \otimes X_{i+n})  \otimes R \otimes \rho^U(X_{i+n+1}, \ldots, X_{m+n-1})\big) \circ\\
& \quad  \big(U \otimes F(X_1) \otimes \ldots \otimes F(X_i) \otimes g_{X_{i+1}, \ldots, X_{i+n}} \otimes U \otimes F(X_{i+n+1})\otimes \ldots \otimes  F(X_{m+n-1}) \big)  \circ \\
& \quad \big(U \otimes F(X_1) \otimes \ldots \otimes F(X_i) \otimes U \otimes \rho^U(X_{i+1}, \ldots, X_{i+n}) \otimes F(X_{i+n+1}) \otimes \ldots \otimes F(X_{m+n-1})\big)\circ \\
& \quad \big(U \otimes F(X_1) \otimes \ldots \otimes F(X_i) \otimes \Delta \otimes  F(X_{i+1}) \otimes \ldots \otimes F(X_{i+n})  \otimes \ldots  F(X_{m+n-1})\big)\circ \\
&\quad \big(U \otimes \rho^U(X_1,\ldots, X_i) \otimes F(X_{i+1})  \otimes \ldots  \otimes F(X_{m+n-1})  \big) \circ (\Delta \otimes F(X_1) \otimes \ldots \otimes F(X_{m+n-1})) 
\end{align*}
Applying the condition in \eqref{3.10cx} for $g\in P^n(F,U,R,\,\overline{\!\Psi\!}\,)$, the above equates to
\begin{align*}
&\big( F(X_1) \otimes \ldots \otimes F(X_{m+n-1}) \otimes \mu \otimes U\big) \circ \big(f_{X_1, \ldots, X_i, X_{i+1} \otimes \ldots \otimes X_{i+n}, \ldots, X_{m+n-1}} \otimes R \otimes U\big) \circ \\
&\big(U \otimes  F(X_1) \otimes \ldots \otimes F(X_i) \otimes F(X_{i+1} \otimes \ldots \otimes X_{i+n}) \otimes \rho^R(X_{i+n+1}, \ldots, X_{m+n-1}) \otimes U\big) \circ \\
& \big(U \otimes  F(X_1) \otimes \ldots \otimes F(X_i) \otimes F(X_{i+1} \otimes \ldots \otimes X_{i+n}) \otimes R \otimes \rho^U(X_{i+n+1}, \ldots, X_{m+n-1})\big)\circ \\
& \big(U \otimes  F(X_1) \otimes \ldots \otimes F(X_i) \otimes F(X_{i+1} \otimes \ldots \otimes X_{i+n}) \otimes \,\overline{\!\Psi\!}\,\otimes F(X_{i+n+1}) \otimes \ldots \otimes F(X_{m+n-1})\big) \circ \\
& \big(U \otimes  F(X_1) \otimes \ldots \otimes F(X_i) \otimes \rho^U(X_{i+1}, \ldots,  X_{i+n}) \otimes R \otimes F(X_{i+n+1}) \otimes \ldots \otimes F(X_{m+n-1})\big) \circ\\
& \big(U \otimes  F(X_1) \otimes \ldots \otimes F(X_i) \otimes U \otimes g_{X_{i+1}, \ldots, X_{i+n}} \otimes F(X_{i+n+1}) \otimes \ldots \otimes F(X_{m+n-1})\big) \circ\\
&   \big(U \otimes  F(X_1) \otimes \ldots \otimes F(X_i) \otimes \Delta \otimes F(X_{i+1}) \otimes \ldots \otimes F(X_{m+n-1}) \big)\circ \\
& \big(U \otimes \rho^U(X_1,\ldots, X_i) \otimes F(X_{i+1})  \otimes \ldots  F(X_{m+n-1})  \big) \circ (\Delta \otimes F(X_1) \otimes \ldots \otimes F(X_{m+n-1})) 
 \end{align*}
which now using \eqref{morp}, becomes
\begin{align*}
&\big( F(X_1) \otimes \ldots \otimes F(X_{m+n-1}) \otimes \mu \otimes U\big) \circ \big(f_{X_1, \ldots, X_i, X_{i+1} \otimes \ldots \otimes X_{i+n}, \ldots, X_{m+n-1}} \otimes R \otimes U\big) \circ\\
&\big(U \otimes  F(X_1) \otimes \ldots \otimes F(X_i) \otimes F(X_{i+1} \otimes \ldots \otimes X_{i+n}) \otimes \ldots \otimes F(X_{m+n-1}) \otimes \overline{\Psi}\big) \circ\\
&\big(U \otimes  F(X_1) \otimes \ldots \otimes F(X_i) \otimes \rho^U(X_{i+1}, \ldots,  X_{m+n-1}) \otimes R\big)\\
&\big(U \otimes  F(X_1) \otimes \ldots \otimes F(X_i) \otimes U \otimes  F(X_{i+1} \otimes \ldots \otimes X_{i+n}) \otimes \rho^R(X_{i+n+1}, \ldots, X_{m+n-1})\big) \circ \\
&\big(U \otimes  F(X_1) \otimes \ldots \otimes F(X_i) \otimes U \otimes g_{X_{i+1},\ldots,X_{i+n}} \otimes F(X_{i+n+1}) \otimes \ldots \otimes F(X_{m+n-1})\big) \circ\\
&\big(U \otimes  F(X_1) \otimes \ldots \otimes F(X_i) \otimes \Delta \otimes F(X_{i+1}) \otimes \ldots \otimes F(X_{m+n-1}) \big)\\
& \big(U \otimes \rho^U(X_1,\ldots, X_i) \otimes F(X_{i+1})  \otimes \ldots  F(X_{m+n-1})  \big) \circ (\Delta \otimes F(X_1) \otimes \ldots \otimes F(X_{m+n-1})) 
\end{align*}
Using \eqref{2.8bv}, the above equates to
\begin{align*}
&\quad \big( F(X_1) \otimes \ldots \otimes F(X_{m+n-1}) \otimes \mu \otimes U\big) \circ \big(f_{X_1, \ldots, X_i, X_{i+1} \otimes \ldots \otimes X_{i+n}, \ldots, X_{m+n-1}} \otimes R \otimes U\big) \circ\\
& \quad \big(U \otimes  F(X_1) \otimes \ldots \otimes F(X_i) \otimes F(X_{i+1} \otimes \ldots \otimes X_{i+n}) \otimes \ldots \otimes F(X_{m+n-1}) \otimes \,\overline{\!\Psi\!}\,\big) \circ \\
& \quad \big(U \otimes  F(X_1) \otimes \ldots \otimes F(X_i) \otimes \rho^U(X_{i+1}, \ldots,  X_{m+n-1}) \otimes R\big) \\
& \quad \big(U \otimes  F(X_1) \otimes \ldots \otimes F(X_i) \otimes U \otimes  F(X_{i+1} \otimes \ldots \otimes X_{i+n}) \otimes \rho^R(X_{i+n+1}, \ldots, X_{m+n-1})\big) \circ \\
& \quad \big(U \otimes  F(X_1) \otimes \ldots \otimes F(X_i) \otimes U \otimes g_{X_{i+1},\ldots,X_{i+n}} \otimes F(X_{i+n+1}) \otimes \ldots \otimes F(X_{m+n-1})\big) \circ\\
& \quad \big(U \otimes \rho^U(X_1,\ldots, X_i) \otimes U \otimes  F(X_{i+1}) \otimes \ldots \otimes F(X_{m+n-1}) \big)\\
&\quad \big(U \otimes U \otimes \rho^U(X_1,\ldots, X_i)  \otimes  F(X_{i+1}) \otimes \ldots \otimes F(X_{m+n-1}) \big)\\
&\quad \big(U \otimes \Delta \otimes F(X_1) \otimes \ldots \otimes F(X_{m+n-1})\big) \circ (\Delta \otimes F(X_1) \otimes \ldots \otimes F(X_{m+n-1})) \\
&=\big( F(X_1) \otimes \ldots \otimes F(X_{m+n-1}) \otimes \mu \otimes U\big) \circ \big(F(X_1) \otimes \ldots \otimes F(X_{m+n-1}) \otimes R \otimes \,\overline{\!\Psi\!}\,\big) \circ \big(f_{X_1, \ldots, X_i, X_{i+1} \otimes \ldots \otimes X_{i+n}, \ldots, X_{m+n-1}} \otimes U \otimes R\big) \circ \\
& \quad \big(U \otimes \rho^U(X_{1}, \ldots,  X_{m+n-1}) \otimes R \big) \circ (\Delta \otimes F(X_1) \otimes \ldots \otimes F(X_{m+n-1}) \otimes R) \circ\\
& \quad \big(U \otimes  F(X_1) \otimes \ldots \otimes F(X_i) \otimes F(X_{i+1})\otimes ...\otimes F(X_{i+n})\otimes \rho^R(X_{i+n+1}, \ldots, X_{m+n-1})\big) \circ\\
& \quad \big(U \otimes  F(X_1) \otimes \ldots \otimes F(X_i) \otimes  g_{X_{i+1},\ldots,X_{i+n}} \otimes F(X_{i+n+1}) \otimes \ldots \otimes F(X_{m+n-1})\big) \circ\\
&\quad \big(U  \otimes \rho^U(X_1,\ldots, X_i)  \otimes  F(X_{i+1}) \otimes \ldots \otimes F(X_{m+n-1}) \big) \circ (\Delta \otimes F(X_1) \otimes \ldots \otimes F(X_{m+n-1})) 
\end{align*}
Again using the commutative diagram in \eqref{3.10cx} for $f$, the above equates to
\begin{align*}
&\quad \big( F(X_1) \otimes \ldots \otimes F(X_{m+n-1}) \otimes \mu \otimes U\big) \circ \big(F(X_1) \otimes \ldots \otimes F(X_{m+n-1}) \otimes R \otimes \,\overline{\!\Psi\!}\,\big) \circ \\
&\quad \big(F(X_1) \otimes \ldots \otimes F(X_i) \otimes  F(X_{i+1} \otimes \ldots \otimes X_{i+n}) \otimes  F(X_{i+n+1}) \otimes \ldots \otimes  F(X_{m+n-1}) \otimes \,\overline{\!\Psi\!}\,\otimes R\big) \circ \big(\rho^U(X_1,\ldots, X_{m+n-1})  \otimes R \otimes R \big) \circ\\
& \quad \big(U \otimes f_{X_1, \ldots, X_i, X_{i+1} \otimes \ldots \otimes X_{i+n}, \ldots, X_{m+n-1}} \otimes R \big) \circ (\Delta \otimes F(X_1) \otimes \ldots \otimes F(X_{m+n-1})\otimes R) \circ 
 \\  &\quad \big(U \otimes F(X_1) \otimes \ldots \otimes F(X_{i+n}) \otimes \rho^R(X_{i+n+1}, \ldots, X_{m+n-1})\big) \circ \\
& \quad \big(U \otimes  F(X_1) \otimes \ldots \otimes F(X_i)  \otimes  g_{X_{i+1},\ldots,X_{i+n}} \otimes F(X_{i+n+1}) \otimes \ldots \otimes F(X_{m+n-1})\big) \circ\\
&\quad \big(U  \otimes \rho^U(X_1,\ldots, X_i)  \otimes  F(X_{i+1}) \otimes \ldots \otimes F(X_{m+n-1}) \big) \circ (\Delta \otimes F(X_1) \otimes \ldots \otimes F(X_{m+n-1})) 
\end{align*}
The result now follows by using the condition in \eqref{mult}.
\end{proof}


We now  construct another subcomplex of $C^\bullet_{DY}(F,U,R)$. For this, we consider a morphism $\,\overline{\!\Phi\!}\,: U \otimes R \longrightarrow R \otimes U$ in the category $\mathcal Z(F)$ which is compatible with the comultiplication $\Delta$ and counit $\epsilon$ of the coalgebra object $U$. This means that in addition to the condition as in \eqref{morp}, the morphism $\overline{\!\Phi\!}\,$ fits into the following commutative diagram
\begin{equation}\label{comult}
\begin{tikzcd}
U \otimes R
\arrow[r, "{\Delta \otimes R}"]
\arrow[d, swap, "{\overline{\!\Phi\!}\,}"]
&
U \otimes U \otimes R
\arrow[r, "{U \otimes \overline{\!\Phi\!}\,}"]
&
U \otimes R \otimes U
\arrow[d, "{\overline{\!\Phi\!}\, \otimes U}"]\\
R \otimes U
\arrow[rr,"{R \otimes \Delta }"]
& 
&
R \otimes U \otimes U
\end{tikzcd}
\qquad \qquad
\begin{tikzcd}
U \otimes R
\arrow[rr,"{\overline{\!\Phi\!}\,}"]
\arrow[dr, swap, "{\epsilon \otimes R }"]
&
&
R \otimes U
\arrow[dl, "{ R \otimes \epsilon}"]\\
& R
\end{tikzcd}
\end{equation}

For each $m \geq 0$, we now set $Q^m(F,U,R,\overline{\!\Phi\!}\,)\subseteq C^m_{DY}(F,U,R)=Nat(U\otimes F^{\otimes m}, F^{\otimes m}\otimes R)$ to be the collection of those natural tranformations $f:U \otimes F^{\otimes m} \longrightarrow F^{\otimes m} \otimes R$ for which the following diagram commutes:

\begin{equation}\label{4.5yy}\small
\begin{tikzcd}[column sep=2.44cm, row sep=large]
U \otimes R \otimes F(X_1) \otimes \ldots \otimes F(X_m)
\arrow[r, "{U \otimes \rho^R(X_1,\ldots,X_m)}"]
\arrow[d, "{\overline{\!\Phi\!}\, \otimes F(X_1) \otimes \ldots \otimes F(X_m)}"']
&
U \otimes F(X_1) \otimes \ldots \otimes F(X_m) \otimes R
\arrow[dd, "{f_{X_1,\ldots,X_m} \otimes R}"]
&
\\
R \otimes U \otimes F(X_1) \otimes \ldots \otimes F(X_m)
\arrow[d, swap, "{R \otimes f_{X_1,\ldots,X_m}}"]
\\
R \otimes F(X_1) \otimes \ldots \otimes F(X_m) \otimes R
\arrow[r, "{\rho^R(X_1,\ldots,X_m) \otimes R}"]
&
F(X_1) \otimes \ldots \otimes F(X_m) \otimes R \otimes R
\arrow[r, "{F(X_1) \otimes \ldots \otimes F(X_m) \otimes \mu}"]
&
F(X_1) \otimes \ldots \otimes F(X_m) \otimes R 
 \end{tikzcd}
 \end{equation}

\begin{prop}\label{prop4.3}
The subspaces $Q^m(F,U,R,\overline{\!\Phi\!}\,)\subseteq C^m_{DY}(F,U,R)$ for $m\geq 0$ determine a subcomplex of $(C^\bullet_{DY}(F,U,R),\partial^\bullet)$.
\end{prop}
\begin{proof} 
This follows in a manner dual to Proposition \ref{ent11}, by applying the condition in \eqref{morp} for the map $~\overline{\!\Phi\!}\,$ and using \eqref{Rcom} and \eqref{4.5yy}.
\end{proof}

The next result says that the operations $\diamond_i$ also restrict to the subcomplex $Q^\bullet(F,U,R,\overline{\!\Phi\!}\,)$.

\begin{lem}\label{closedn}
For $f \in Q^m(F,U,R,\overline{\!\Phi\!}\,)$ and $g \in Q^n(F,U,R,\overline{\!\Phi\!}\,)$, we have each $f \diamond_i g \in Q^{m+n-1}(F,U,R,\overline{\!\Phi\!}\,)$.
\end{lem}
\begin{proof}
This follows  in a manner dual to Lemma \ref{closed},  using the associativity of $\mu$, the condition in \eqref{4.5yy} and the fact that  $~\overline{\!\Phi\!}\,$ is a morphism in $\mathcal{Z}(F)$.
\end{proof}

We are now ready to construct  a subcomplex of  $\big(C^\bullet_{DY}(F,U,R), \partial^\bullet\big)$ which corresponds to a comp algebra, and whose cohomology therefore carries the structure of a Gerstenhaber algebra. For this, we consider an entwining structure $~\overline{\!\Theta\!}\,:U \otimes R \longrightarrow R \otimes U$ in the monoidal category $\mathcal{Z}(F)$, which means that $~\overline{\!\Theta\!}\,$ is a morphism in $\mathcal{Z}(F)$ satisfying both the conditions in \eqref{mult} and in \eqref{comult}. In other words, the morphism $~\overline{\!\Theta\!}\,$ in $\mathcal{Z}(F)$  is compatible with the comultiplication $\Delta$ and counit $\epsilon$ on the coalgebra object $(U,\rho^U,\Delta,\epsilon)$ as well as the multiplication $\mu$ and the unit $\eta$ on the algebra object $(R, \rho^R,\mu,\eta)$. We   consider the corresponding subcomplexes $\big(P^\bullet(F,U,R,~\overline{\!\Theta\!}\,), \partial^\bullet\big)$ and $\big(Q^\bullet(F,U,R,~\overline{\!\Theta\!}\,),\partial^\bullet\big)$ of $\big(C^\bullet_{DY}(F,U,R), \partial^\bullet\big)$ and set, for each $m \geq 0$,
\begin{equation}
D^m(F,U,R,~\overline{\!\Theta\!}\,):= P^m(F,U,R,~\overline{\!\Theta\!}\,) \cap Q^m(F,U,R,~\overline{\!\Theta\!}\,)\subseteq C^m_{DY}(F,U,R)
\end{equation}
 It is   clear that $\big(D^\bullet(F,U,R,~\overline{\!\Theta\!}\,), \partial^\bullet\big)$ is a subcomplex of $\big(C^\bullet_{DY}(F,U,R), \partial^\bullet\big)$. We now have the following result.
\begin{Thm}\label{rysz}
 Let $F:\mathscr C\longrightarrow \mathscr D$ be a monoidal  functor between monoidal  categories, $(U,\rho^U,\Delta,\epsilon)$ be a coalgebra object and $(R, \rho^R,\mu,\eta)$ be an algebra object in the centralizer $\mathcal Z(F)$ of $F$. Let $~\overline{\!\Theta\!}\,: U \otimes R \longrightarrow R \otimes U$ be an entwining structure in $\mathcal{Z}(F)$. Then,

\smallskip 
(a)  $P^\bullet(F,U,R,~\overline{\!\Theta\!}\,)$ and $Q^\bullet(F,U,R,~\overline{\!\Theta\!}\,)$ are weak comp algebras.

\smallskip
(b)  $D^\bullet(F,U,R,~\overline{\!\Theta\!}\,)$ is a comp algebra.
\end{Thm}
\begin{proof}
(a)  From  Lemma \ref{closed}, we know that 
$P^\bullet(F,U,R, ~\overline{\!\Theta\!}\,)$ is closed under the operations $\diamond_i$. We will now check that the element $\pi \in C^2_{DY}(F,U,R)$ defined in \eqref{pii} actually lies in  $P^2(F,U,R,~\overline{\!\Theta\!}\,)$. For this, we note that if $X_1$, $X_2\in Ob(\mathscr C)$, we have
\begin{align*}
&(\pi_{X_1,X_2} \otimes U) \circ (U \otimes \rho^U(X_1,X_2)) \circ (\Delta \otimes F(X_1) \otimes F(X_2))\\
&=(F(X_1) \otimes F(X_2) \otimes \eta \otimes U) \circ (\epsilon \otimes F(X_1) \otimes F(X_2) \otimes U) \circ (U \otimes \rho^U(X_1,X_2)) \circ (\Delta \otimes F(X_1) \otimes F(X_2))\\
&=(F(X_1) \otimes F(X_2) \otimes \eta \otimes U) \circ \rho^U(X_1,X_2) \circ (\epsilon \otimes U \otimes F(X_1) \otimes F(X_2)) \circ (\Delta \otimes F(X_1) \otimes F(X_2))\\
&=(F(X_1) \otimes F(X_2) \otimes ~\overline{\!\Theta\!}\,) \circ (F(X_1) \otimes F(X_2) \otimes U \otimes \eta) \circ \rho^U(X_1,X_2) \circ (U \otimes \epsilon \otimes F(X_1) \otimes F(X_2)) \circ (\Delta \otimes F(X_1) \otimes F(X_2))\\
&=(F(X_1) \otimes F(X_2) \otimes ~\overline{\!\Theta\!}\,) \circ (\rho^U(X_1,X_2) \otimes R) \circ ( U \otimes F(X_1) \otimes F(X_2) \otimes \eta) \circ (U \otimes \epsilon \otimes F(X_1) \otimes F(X_2)) \circ (\Delta \otimes F(X_1) \otimes F(X_2))\\
&=(F(X_1) \otimes F(X_2) \otimes ~\overline{\!\Theta\!}\,) \circ (\rho^U(X_1,X_2) \otimes R) \circ (U \otimes \pi_{X_1,X_2}) \circ (\Delta \otimes F(X_1) \otimes F(X_2))
\end{align*}
This shows that $\pi \in C^2_{DY}(F,U,R)$ satisfies the condition in \eqref{3.10cx}, i.e.,  $\pi \in P^{2}(F,U,R,~\overline{\!\Theta\!}\,)$. Since $C^\bullet_{DY}(F,U,R)$ is a weak comp algebra by Proposition \ref{p3.4f}, it is now clear that $ P^{\bullet}(F,U,R,~\overline{\!\Theta\!}\,) \subseteq C^\bullet_{DY}(F,U,R)$ is   a weak comp algebra. Similarly, it may be verified that $\pi \in Q^{2}(F,U,R,~\overline{\!\Theta\!}\,)$. Combining with Lemma \ref{closedn},  it  follows that $ Q^{\bullet}(F,U,R,~\overline{\!\Theta\!}\,) \subseteq C^\bullet_{DY}(F,U,R)$ is   also a weak comp algebra.


\smallskip
(b)  Using Lemma \ref{closed} and Lemma \ref{closedn}, it is clear that the intersection 
$D^\bullet(F,U,R, ~\overline{\!\Theta\!}\,)$ is closed under the operations $\diamond_i$. From the proof of part (a), it follows that $\pi \in D^{2}(F,U,R,~\overline{\!\Theta\!}\,)=P^{2}(F,U,R,~\overline{\!\Theta\!}\,)
\cap Q^{2}(F,U,R,~\overline{\!\Theta\!}\,)$. Hence, $D^\bullet(F,U,R,~\overline{\!\Theta\!}\,)$ is a weak comp algebra. To show that $D^\bullet(F,U,R,~\overline{\!\Theta\!}\,)$ is a comp algebra, it remains to check that the condition 
\begin{equation} (f \diamond_i g) \diamond_j h=(f \diamond_j h) \diamond_{i+p-1} g
\end{equation} in Definition \ref{wg} is satisfied for all $f \in D^m(F,U,R, ~\overline{\!\Theta\!}\,)$, $g \in D^n(F,U,R,~\overline{\!\Theta\!}\,)$ and $h \in D^p(F,U,R,~\overline{\!\Theta\!}\,)$ whenever $j<i$. Let $X_1, \ldots, X_{m+n+p-2} \in Ob(\mathscr C)$. First using \eqref{2.8bv}, we see that
\begin{align*}\label{fbl4}
&((f \diamond_j h) \diamond_{i+p-1} g)_{X_1,...,X_{m+n+p-2}}\\
&=\big( F(X_1) \otimes \ldots \otimes F(X_{m+n+p-2}) \otimes \mu\big) \circ \big( F(X_1) \otimes \ldots \otimes F(X_{i+p-1}) \otimes F(X_{i+p} \otimes \ldots \otimes X_{i+p+n-1}) \otimes F(X_{i+p+n}) \otimes \ldots \otimes F(X_{m+n+p-2}) \otimes \mu \otimes R\big) \\
& \quad \circ \big(f_{X_1,\ldots, X_{j+1}\otimes \ldots \otimes X_{j+p}, \ldots, X_{p+i-1}, X_{p+i} \otimes \ldots \otimes X_{n+p+i-1}, \ldots, X_{m+n+p-2}} \otimes R \otimes R\big)\circ \\
& \quad\big(U \otimes F(X_1) \otimes \ldots \otimes F(X_{j}) \otimes F(X_{j+1} \otimes \ldots \otimes X_{j+p}) \otimes \rho^R(X_{j+p+1}, \ldots, X_{i+p} \otimes \ldots \otimes X_{i+p+n-1}, \ldots, X_{m+n+p-2}) \otimes R\big)\\
& \quad \big(U \otimes F(X_1) \otimes \ldots \otimes F(X_j) \otimes h_{X_{j+1}, \ldots, X_{j+p}} \otimes F(X_{p+j+1}) \otimes \ldots \otimes F(X_{i+p} \otimes \ldots \otimes X_{i+p+n-1}) \otimes \ldots \otimes F(X_{m+n+p-2}) \otimes R\big)  \circ \\
& \quad \big(U \otimes F(X_1) \otimes \ldots \otimes F(X_j) \otimes U \otimes F(X_{j+1}) \otimes \ldots \otimes F(X_{i+p+n-1}) \otimes \rho^R(X_{i+p+n}, \ldots, X_{m+n+p-2})\big) \circ\\
& \quad \big(U \otimes F(X_1) \otimes \ldots \otimes F(X_j) \otimes U \otimes F(X_{j+1}) \otimes \ldots \otimes F(X_{i+p-1}) \otimes g_{X_{i+p}, \ldots, X_{i+p+n-1}} \otimes F(X_{i+p+n})\otimes \ldots \otimes  F(X_{m+n+p-2}) \big)  \circ \\
& \quad \big(U \otimes F(X_1) \otimes \ldots \otimes F(X_j) \otimes U \otimes \rho^U(X_{j+1}, \ldots, X_{i+p-1}) \otimes F(X_{i+p}) \otimes \ldots \otimes F(X_{m+n+p-2})\big)  \circ \\
& \quad \big(U \otimes F(X_1) \otimes \ldots \otimes F(X_j) \otimes \Delta \otimes F(X_{j+1}) \otimes \ldots \otimes F(X_{m+n+p-2})\big)  \circ \big(U \otimes \rho^U(X_1, \ldots, X_j) \otimes F(X_{j+1}) \otimes \ldots \otimes F(X_{m+n+p-2})\big) \circ \\
& \quad (\Delta \otimes F(X_1) \otimes \ldots \otimes F(X_{m+n+p-2})) \\ 
&=\big( F(X_1) \otimes \ldots \otimes F(X_{m+n+p-2}) \otimes \mu\big) \circ \big( F(X_1) \otimes \ldots \otimes F(X_{i+p-1}) \otimes F(X_{i+p} \otimes \ldots \otimes X_{i+p+n-1}) \otimes F(X_{i+p+n}) \otimes \ldots \otimes F(X_{m+n+p-2}) \otimes \mu \otimes R\big) \\
& \quad \circ \big(f_{X_1,\ldots, X_{j+1}\otimes \ldots \otimes X_{j+p}, \ldots, X_{p+i-1}, X_{p+i} \otimes \ldots \otimes X_{n+p+i-1}, \ldots, X_{m+n+p-2}} \otimes R \otimes R\big)\circ \\
& \quad\big(U \otimes F(X_1) \otimes \ldots \otimes F(X_{j}) \otimes F(X_{j+1} \otimes \ldots \otimes X_{j+p}) \otimes \rho^R(X_{j+p+1}, \ldots, X_{i+p} \otimes \ldots \otimes X_{i+p+n-1}, \ldots, X_{m+n+p-2}) \otimes R\big)\\
& \quad\big(U \otimes F(X_1) \otimes \ldots \otimes F(X_{j+p}) \otimes R \otimes F(X_{j+p+1}) \otimes \ldots \otimes F(X_{i+p-1}) \otimes F(X_{i+p} \otimes \ldots \otimes X_{i+p+n-1}) \otimes \rho^R(X_{i+p+n}, \ldots, X_{m+n+p-2})\big) \circ\\
& \quad \big(U \otimes F(X_1) \otimes \ldots \otimes F(X_{j+p}) \otimes R \otimes F(X_{j+p+1}) \otimes \ldots \otimes  F(X_{i+p-1}) \otimes g_{X_{i+p}, \ldots, X_{i+p+n-1}} \otimes F(X_{i+p+n})\otimes \ldots \otimes  F(X_{m+n+p-2}) \big)  \circ \\
& \quad \big(U \otimes F(X_1) \otimes \ldots \otimes F(X_{j+p}) \otimes R \otimes \rho^U(X_{j+p+1}, \ldots, X_{i+p-1}) \otimes F(X_{i+p})\otimes \ldots \otimes  F(X_{m+n+p-2}) \big)  \circ \\
& \quad \big(U \otimes F(X_1) \otimes \ldots \otimes F(X_{j}) \otimes h_{X_{j+1}, \ldots, X_{j+p}} \otimes U \otimes F(X_{p+j+1}) \otimes \ldots \otimes F(X_{m+n+p-2})\big)  \circ \\
& \quad \big(U \otimes F(X_1) \otimes \ldots \otimes F(X_{j}) \otimes U \otimes  \rho^U(X_{j+1}, \ldots, X_{j+p}) \otimes F(X_{p+j+1}) \otimes \ldots \otimes F(X_{m+n+p-2})\big)  \circ \\
& \quad \big(U \otimes F(X_1) \otimes \ldots \otimes F(X_{j}) \otimes \Delta \otimes  F(X_{j+1}) \otimes \ldots \otimes F(X_{m+n+p-2})\big)  \circ \big(U \otimes \rho^U(X_1, \ldots, X_j) \otimes F(X_{j+1}) \otimes \ldots \otimes F(X_{m+n+p-2})\big) \circ \\
& \quad (\Delta \otimes F(X_1) \otimes \ldots \otimes F(X_{m+n+p-2})) 
\end{align*}
Now using the commutative diagram as in \eqref{3.10cx} for the maps $h$ and $~\overline{\!\Theta\!}\,$, the above equates to
\begin{align*}
&\quad \big( F(X_1) \otimes \ldots \otimes F(X_{m+n+p-2}) \otimes \mu\big) \circ \big( F(X_1) \otimes \ldots \otimes F(X_{i+p-1}) \otimes F(X_{i+p} \otimes \ldots \otimes X_{i+p+n-1}) \otimes F(X_{i+p+n}) \otimes \ldots \otimes F(X_{m+n+p-2}) \otimes \mu \otimes R\big) \\
& \quad \circ \big(f_{X_1,\ldots, X_{j+1}\otimes \ldots \otimes X_{j+p}, \ldots, X_{p+i-1}, X_{p+i} \otimes \ldots \otimes X_{n+p+i-1}, \ldots, X_{m+n+p-2}} \otimes R \otimes R\big)\circ \\
& \quad\big(U \otimes F(X_1) \otimes \ldots \otimes F(X_{j}) \otimes F(X_{j+1} \otimes \ldots \otimes X_{j+p}) \otimes \rho^R(X_{j+p+1}, \ldots, X_{i+p} \otimes \ldots \otimes X_{i+p+n-1}, \ldots, X_{m+n+p-2}) \otimes R\big)\\
& \quad\big(U \otimes F(X_1) \otimes \ldots \otimes F(X_{j+p}) \otimes R \otimes F(X_{j+p+1}) \otimes \ldots \otimes F(X_{i+p-1}) \otimes F(X_{i+p} \otimes \ldots \otimes X_{i+p+n-1}) \otimes \rho^R(X_{i+p+n}, \ldots, X_{m+n+p-2})\big) \circ\\
& \quad \big(U \otimes F(X_1) \otimes \ldots \otimes F(X_{j+p}) \otimes R \otimes F(X_{j+p+1}) \otimes \ldots \otimes  F(X_{i+p-1}) \otimes g_{X_{i+p}, \ldots, X_{i+p+n-1}} \otimes F(X_{i+p+n})\otimes \ldots \otimes  F(X_{m+n+p-2}) \big)  \circ \\
& \quad \big(U \otimes F(X_1) \otimes \ldots \otimes F(X_{j+p}) \otimes R \otimes \rho^U(X_{j+p+1}, \ldots, X_{i+p-1}) \otimes F(X_{i+p})\otimes \ldots \otimes  F(X_{m+n+p-2}) \big)  \circ \\
& \quad \big(U \otimes F(X_1) \otimes \ldots \otimes F(X_j) \otimes F(X_{j+1} \otimes \ldots \otimes X_{j+p}) \otimes ~\overline{\!\Theta\!}\, \otimes F(X_{j+p+1})\otimes \ldots \otimes  F(X_{m+n+p-2}) \big)  \circ \\
& \quad \big(U \otimes F(X_1) \otimes \ldots \otimes F(X_{j}) \otimes  \rho^U(X_{j+1}, \ldots, X_{j+p}) \otimes R \otimes F(X_{j+p+1}) \otimes \ldots \otimes F(X_{m+n+p-2})\big)  \circ \\
& \quad \big(U \otimes F(X_1) \otimes \ldots \otimes F(X_{j}) \otimes U \otimes h_{X_{j+1}, \ldots, X_{j+p}}  \otimes F(X_{j+p+1}) \otimes \ldots \otimes F(X_{m+n+p-2})\big)  \circ \\
& \quad \big(U \otimes F(X_1) \otimes \ldots \otimes F(X_{j}) \otimes \Delta \otimes F(X_{j+1}) \otimes \ldots \otimes  F(X_{j+p}) \otimes \ldots \otimes F(X_{m+n+p-2})\big)  \circ \\
& \quad \big(U \otimes \rho^U(X_1, \ldots, X_j) \otimes F(X_{j+1}) \otimes \ldots \otimes F(X_{m+n+p-2})\big) \circ (\Delta \otimes F(X_1) \otimes \ldots \otimes F(X_{m+n+p-2})) 
\end{align*}
which is further equal to
\begin{align*}
&\big( F(X_1) \otimes \ldots \otimes F(X_{m+n+p-2}) \otimes \mu\big) \circ \big( F(X_1) \otimes \ldots \otimes F(X_{i+p-1}) \otimes F(X_{i+p} \otimes \ldots \otimes X_{i+p+n-1}) \otimes F(X_{i+p+n}) \otimes \ldots \otimes F(X_{m+n+p-2}) \otimes \mu \otimes R\big) \\
& \circ \big(f_{X_1,\ldots, X_{j+1}\otimes \ldots \otimes X_{j+p}, \ldots, X_{p+i-1}, X_{p+i} \otimes \ldots \otimes X_{n+p+i-1}, \ldots, X_{m+n+p-2}} \otimes R \otimes R\big)\circ \big(U \otimes F(X_1) \otimes \ldots \otimes F(X_{i+p-1}) \otimes \rho^R(X_{p+i}, \ldots, X_{m+n+p-2}) \otimes R\big) \circ\\
& \big(U \otimes F(X_1) \otimes \ldots \otimes F(X_{i+p-1}) \otimes R \otimes F(X_{i+p} \otimes \ldots \otimes X_{i+p+n-1}) \otimes \rho^R(X_{i+p+n}, \ldots, X_{m+n+p-2})\big) \circ\\
& \big(U \otimes F(X_1) \otimes \ldots \otimes F(X_{i+p-1}) \otimes R \otimes g_{X_{i+p}, \ldots, X_{i+p+n-1}} \otimes F(X_{i+p+n})\otimes \ldots \otimes  F(X_{m+n+p-2}) \big)  \circ \\
& \big(U \otimes F(X_1) \otimes \ldots \otimes F(X_{j+p}) \otimes \rho^R(X_{j+p+1}, \ldots, X_{i+p-1} ) \otimes U \otimes F(X_{i+p})\otimes \ldots \otimes  F(X_{m+n+p-2}) \big)  \circ \\
&  \big(U \otimes F(X_1) \otimes \ldots \otimes F(X_{j+p}) \otimes R \otimes \rho^U(X_{j+p+1}, \ldots, X_{i+p-1}) \otimes F(X_{i+p})\otimes \ldots \otimes  F(X_{m+n+p-2}) \big)  \circ \\
&  \big(U \otimes F(X_1) \otimes \ldots \otimes F(X_j) \otimes F(X_{j+1} \otimes \ldots \otimes X_{j+p}) \otimes ~\overline{\!\Theta\!}\, \otimes F(X_{j+p+1})\otimes \ldots \otimes  F(X_{m+n+p-2}) \big)  \circ \\
& \big(U \otimes F(X_1) \otimes \ldots \otimes F(X_{j}) \otimes  \rho^U(X_{j+1}, \ldots, X_{j+p}) \otimes R \otimes F(X_{j+p+1}) \otimes \ldots \otimes F(X_{m+n+p-2})\big)  \circ \\
&  \big(U \otimes F(X_1) \otimes \ldots \otimes F(X_{j}) \otimes U \otimes h_{X_{j+1}, \ldots, X_{j+p}}  \otimes F(X_{j+p+1}) \otimes \ldots \otimes F(X_{m+n+p-2})\big)  \circ \\
&  \big(U \otimes F(X_1) \otimes \ldots \otimes F(X_{j}) \otimes \Delta \otimes F(X_{j+1}) \otimes \ldots \otimes  F(X_{j+p}) \otimes \ldots \otimes F(X_{m+n+p-2})\big)  \circ \\
& \big(U \otimes \rho^U(X_1, \ldots, X_j) \otimes F(X_{j+1}) \otimes \ldots \otimes F(X_{m+n+p-2})\big) \circ (\Delta \otimes F(X_1) \otimes \ldots \otimes F(X_{m+n+p-2})) 
\end{align*}
Using the fact that $~\overline{\!\Theta\!}\,$ is a morphism in $\mathcal{Z}(F)$, the above becomes
\begin{align*}
&\big( F(X_1) \otimes \ldots \otimes F(X_{m+n+p-2}) \otimes \mu\big) \circ \big( F(X_1) \otimes \ldots \otimes F(X_{i+p-1}) \otimes F(X_{i+p} \otimes \ldots \otimes X_{i+p+n-1}) \otimes F(X_{i+p+n}) \otimes \ldots \otimes F(X_{m+n+p-2}) \otimes \mu \otimes R\big) \\
& \circ \big(f_{X_1,\ldots, X_{j+1}\otimes \ldots \otimes X_{j+p}, \ldots, X_{p+i-1}, X_{p+i} \otimes \ldots \otimes X_{n+p+i-1}, \ldots, X_{m+n+p-2}} \otimes R \otimes R\big)\circ \big(U \otimes F(X_1) \otimes \ldots \otimes F(X_{i+p-1}) \otimes \rho^R(X_{p+i}, \ldots, X_{m+n+p-2}) \otimes R\big) \circ\\
&  \big(U \otimes F(X_1) \otimes \ldots \otimes F(X_{i+p-1}) \otimes R \otimes F(X_{i+p} \otimes \ldots \otimes X_{i+p+n-1}) \otimes \rho^R(X_{i+p+n}, \ldots, X_{m+n+p-2})\big) \circ\\
& \big(U \otimes F(X_1) \otimes \ldots \otimes F(X_{i+p-1}) \otimes R \otimes g_{X_{i+p}, \ldots, X_{i+p+n-1}} \otimes F(X_{i+p+n})\otimes \ldots \otimes  F(X_{m+n+p-2}) \big)  \circ \\
&  \big(U \otimes F(X_1) \otimes \ldots \otimes F(X_{i+p-1}) \otimes ~\overline{\!\Theta\!}\, \otimes F(X_{i+p}) \otimes \ldots \otimes  F(X_{m+n+p-2}) \big)  \circ \\
& \big(U \otimes F(X_1) \otimes \ldots \otimes F(X_{j+p}) \otimes \rho^U(X_{j+p+1}, \ldots, X_{i+p-1}) \otimes R \otimes F(X_{i+p})\otimes \ldots \otimes  F(X_{m+n+p-2}) \big)  \circ \\
& \big(U \otimes F(X_1) \otimes \ldots \otimes F(X_{j+p}) \otimes U \otimes \rho^R(X_{j+p+1}, \ldots, X_{i+p-1})  \otimes F(X_{i+p}) \otimes \ldots \otimes  F(X_{m+n+p-2}) \big)  \circ \\
&  \big(U \otimes F(X_1) \otimes \ldots \otimes F(X_{j}) \otimes  \rho^U(X_{j+1}, \ldots, X_{j+p}) \otimes R \otimes F(X_{j+p+1}) \otimes \ldots \otimes F(X_{m+n+p-2})\big)  \circ \\
& \big(U \otimes F(X_1) \otimes \ldots \otimes F(X_{j}) \otimes U \otimes h_{X_{j+1}, \ldots, X_{j+p}}  \otimes F(X_{j+p+1}) \otimes \ldots \otimes F(X_{m+n+p-2})\big)  \circ \\
& \big(U \otimes F(X_1) \otimes \ldots \otimes F(X_{j}) \otimes \Delta \otimes F(X_{j+1}) \otimes \ldots \otimes  F(X_{j+p}) \otimes \ldots \otimes F(X_{m+n+p-2})\big)  \circ \\
&  \big(U \otimes \rho^U(X_1, \ldots, X_j) \otimes F(X_{j+1}) \otimes \ldots \otimes F(X_{m+n+p-2})\big) \circ (\Delta \otimes F(X_1) \otimes \ldots \otimes F(X_{m+n+p-2})) 
\end{align*}
which further reduces to
\begin{align*}
& \big( F(X_1) \otimes \ldots \otimes F(X_{m+n+p-2}) \otimes \mu\big) \circ \big(f_{X_1,\ldots, X_{j+1}\otimes \ldots \otimes X_{j+p}, \ldots, X_{p+i-1}, X_{p+i} \otimes \ldots \otimes X_{n+p+i-1}, \ldots, X_{m+n+p-2}} \otimes R \big)\circ\\
& \big(U \otimes F(X_1) \otimes \ldots \otimes F(X_{i+p+n-1}) \otimes \rho^R(X_{i+p+n}, \ldots, X_{m+n+p-2})\big) \circ \\
&\big(U \otimes F(X_1) \otimes \ldots \otimes F(X_{i+p-1}) \otimes F(X_{p+i} \otimes \ldots \otimes X_{i+p+n-1}) \otimes \mu \otimes F(X_{i+p+n}) \otimes \ldots \otimes F(X_{m+n+p-2}) \big) \circ\\
&\big(U \otimes F(X_1) \otimes \ldots \otimes F(X_{i+p-1}) \otimes \rho^R(X_{i+p} \otimes \ldots \otimes X_{i+p+n-1}) \otimes R \otimes F(X_{i+p+n}) \otimes \ldots \otimes F(X_{m+n+p-2}) \big) \circ\\
&\big(U \otimes F(X_1) \otimes \ldots \otimes F(X_{i+p-1}) \otimes R \otimes g_{X_{i+p}, \ldots, X_{i+p+n-1}} \otimes F(X_{i+p+n})\otimes \ldots \otimes  F(X_{m+n+p-2}) \big)  \circ \\
& \big(U \otimes F(X_1) \otimes \ldots \otimes F(X_{i+p-1}) \otimes ~\overline{\!\Theta\!}\, \otimes  F(X_{p+i}) \otimes \ldots \otimes F(X_{m+n+p-2})\big) \circ\\
&\big(U \otimes F(X_1) \otimes \ldots \otimes F(X_{j+p}) \otimes \rho^U(X_{j+p+1}, \ldots, X_{i+p-1}) \otimes R \otimes F(X_{i+p})\otimes \ldots \otimes  F(X_{m+n+p-2}) \big)  \circ \\
& \big(U \otimes F(X_1) \otimes \ldots \otimes F(X_{j+p}) \otimes U \otimes \rho^R(X_{j+p+1}, \ldots, X_{i+p-1})  \otimes F(X_{i+p}) \otimes \ldots \otimes  F(X_{m+n+p-2}) \big)  \circ \\
&\big(U \otimes F(X_1) \otimes \ldots \otimes F(X_{j}) \otimes  \rho^U(X_{j+1}, \ldots, X_{j+p}) \otimes R \otimes F(X_{j+p+1}) \otimes \ldots \otimes F(X_{m+n+p-2})\big)  \circ \\
&\big(U \otimes F(X_1) \otimes \ldots \otimes F(X_{j}) \otimes U \otimes h_{X_{j+1}, \ldots, X_{j+p}}  \otimes F(X_{j+p+1}) \otimes \ldots \otimes F(X_{m+n+p-2})\big)  \circ \\
 &\big(U \otimes F(X_1) \otimes \ldots \otimes F(X_{j}) \otimes \Delta \otimes F(X_{j+1}) \otimes \ldots \otimes  F(X_{j+p}) \otimes \ldots \otimes F(X_{m+n+p-2})\big)  \circ \\
&\big(U \otimes \rho^U(X_1, \ldots, X_j) \otimes F(X_{j+1}) \otimes \ldots \otimes F(X_{m+n+p-2})\big) \circ (\Delta \otimes F(X_1) \otimes \ldots \otimes F(X_{m+n+p-2})) 
\end{align*}
Now we use the condition \eqref{4.5yy} for the map $g$, so that the above becomes
\begin{align*}
&\big( F(X_1) \otimes \ldots \otimes F(X_{m+n+p-2}) \otimes \mu\big) \circ \big(f_{X_1,\ldots, X_{j+1}\otimes \ldots \otimes X_{j+p}, \ldots, X_{p+i-1}, X_{p+i} \otimes \ldots \otimes X_{n+p+i-1}, \ldots, X_{m+n+p-2}} \otimes R \big)\circ \\
&  \big(U \otimes F(X_1) \otimes \ldots \otimes F(X_{i+p+n-1}) \otimes \rho^R(X_{i+p+n}, \ldots, X_{m+n+p-2})\big) \circ \\
& \big(U \otimes F(X_1) \otimes \ldots \otimes F(X_{i+p-1}) \otimes F(X_{p+i} \otimes \ldots \otimes X_{i+p+n-1}) \otimes \mu \otimes F(X_{i+p+n}) \otimes \ldots \otimes F(X_{m+n+p-2}) \big) \circ\\
&  \big(U \otimes F(X_1) \otimes \ldots \otimes F(X_{i+p-1})  \otimes g_{X_{i+p}, \ldots, X_{i+p+n-1}} \otimes R \otimes F(X_{i+p+n})\otimes \ldots \otimes  F(X_{m+n+p-2}) \big)  \circ \\
& \big(U \otimes F(X_1) \otimes \ldots \otimes F(X_{i+p-1}) \otimes U \otimes \rho^R(X_{i+p} \otimes \ldots \otimes X_{i+p+n-1})  \otimes F(X_{i+p+n}) \otimes \ldots \otimes F(X_{m+n+p-2}) \big) \circ\\
& \big(U \otimes F(X_1) \otimes \ldots \otimes F(X_{j+p}) \otimes \rho^U(X_{j+p+1}, \ldots, X_{i+p-1}) \otimes R \otimes F(X_{i+p})\otimes \ldots \otimes  F(X_{m+n+p-2}) \big)  \circ \\
& \big(U \otimes F(X_1) \otimes \ldots \otimes F(X_{j+p}) \otimes U \otimes \rho^R(X_{j+p+1}, \ldots, X_{i+p-1})  \otimes F(X_{i+p}) \otimes \ldots \otimes  F(X_{m+n+p-2}) \big)  \circ \\
& \big(U \otimes F(X_1) \otimes \ldots \otimes F(X_{j}) \otimes  \rho^U(X_{j+1}, \ldots, X_{j+p}) \otimes R \otimes F(X_{j+p+1}) \otimes \ldots \otimes F(X_{m+n+p-2})\big)  \circ \\
&  \big(U \otimes F(X_1) \otimes \ldots \otimes F(X_{j}) \otimes U \otimes h_{X_{j+1}, \ldots, X_{j+p}}  \otimes F(X_{j+p+1}) \otimes \ldots \otimes F(X_{m+n+p-2})\big)  \circ \\
& \big(U \otimes F(X_1) \otimes \ldots \otimes F(X_{j}) \otimes \Delta \otimes F(X_{j+1}) \otimes \ldots \otimes  F(X_{j+p}) \otimes \ldots \otimes F(X_{m+n+p-2})\big)  \circ \\
& \big(U \otimes \rho^U(X_1, \ldots, X_j) \otimes F(X_{j+1}) \otimes \ldots \otimes F(X_{m+n+p-2})\big) \circ (\Delta \otimes F(X_1) \otimes \ldots \otimes F(X_{m+n+p-2})) 
\end{align*}
which further reduces to
\begin{align*}
&\big(F(X_1) \otimes \ldots \otimes F(X_{m+n+p-2}) \otimes \mu\big) \circ \big( F(X_1) \otimes \ldots \otimes F(X_{m+n+p-2}) \otimes R \otimes \mu\big) \circ\\
& \big(f_{X_1,\ldots, X_{j+1}\otimes \ldots \otimes X_{j+p}, \ldots, X_{p+i-1}, X_{p+i} \otimes \ldots \otimes X_{n+p+i-1}, \ldots, X_{m+n+p-2}} \otimes R \otimes R \big)\circ \\
&  \big(U \otimes F(X_1) \otimes \ldots \otimes F(X_{i+p+n-1}) \otimes \rho^R(X_{i+p+n}, \ldots, X_{m+n+p-2}) \otimes R\big) \circ \\
&  \big(U \otimes F(X_1) \otimes \ldots \otimes F(X_{i+p-1}) \otimes g_{X_{i+p}, \ldots, X_{i+p+n-1}} \otimes F(X_{i+p+n})\otimes \ldots \otimes  F(X_{m+n+p-2}) \otimes R \big)  \circ \\
&  \big(U \otimes F(X_1) \otimes \ldots \otimes F(X_{j+p}) \otimes \rho^U(X_{j+p+1}, \ldots, X_{i+p-1}) \otimes F(X_{i+p})\otimes \ldots \otimes  F(X_{m+n+p-2}) \otimes R\big)  \circ \\
&  \big(U \otimes F(X_1) \otimes \ldots \otimes F(X_{j+p}) \otimes U \otimes \rho^R(X_{j+p+1}, \ldots, X_{m+n+p-2})\big)  \circ \\
& \big(U \otimes F(X_1) \otimes \ldots \otimes F(X_{j}) \otimes  \rho^U(X_{j+1}, \ldots, X_{j+p}) \otimes R \otimes F(X_{j+p+1}) \otimes \ldots \otimes F(X_{m+n+p-2})\big)  \circ \\
&  \big(U \otimes F(X_1) \otimes \ldots \otimes F(X_{j}) \otimes U \otimes h_{X_{j+1}, \ldots, X_{j+p}}  \otimes F(X_{j+p+1}) \otimes \ldots \otimes F(X_{m+n+p-2})\big)  \circ \\
&  \big(U \otimes F(X_1) \otimes \ldots \otimes F(X_{j}) \otimes \Delta \otimes F(X_{j+1}) \otimes \ldots \otimes  F(X_{j+p}) \otimes \ldots \otimes F(X_{m+n+p-2})\big)  \circ \\
& \big(U \otimes \rho^U(X_1, \ldots, X_j) \otimes F(X_{j+1}) \otimes \ldots \otimes F(X_{m+n+p-2})\big) \circ (\Delta \otimes F(X_1) \otimes \ldots \otimes F(X_{m+n+p-2}))\\
=~&(f \diamond_i g) \diamond_j h
\end{align*}
where the last equality follows using the compatibility of $~\overline{\!\Theta\!}\,$ and $\mu$ as in \eqref{mult}.
\end{proof}

Let $\diamond:C^m_{DY}(F, U, R) \otimes C^n_{DY}(F, U, R) \longrightarrow C^{m+n-1}_{DY}(F, U, R)$ be the operation defined by setting 
\begin{equation}\label{diamb}
f \diamond g:=\sum\limits_{i=0}^{m-1}  (-1)^{i(n-1)} f \diamond_i g
\end{equation} for $f\in C^m_{DY}(F, U, R) $, $g\in C^n_{DY}(F, U, R) $. From the proof of Theorem \ref{rysz}, we know that $D^\bullet(F,U,R,~\overline{\!\Theta\!}\,)$ is closed under the operations $\diamond_i$ and hence
the operation $\diamond$ as defined in \eqref{diamb} restricts to an operation on $D^\bullet(F,U,R,~\overline{\!\Theta\!}\,)$. We will denote by $H^\bullet(F, U, R, ~\overline{\!\Theta\!}\,)$ the cohomology groups of the subcomplex  $\big(D^\bullet(F,U,R,~\overline{\!\Theta\!}\,),\partial^\bullet\big)$.

\begin{prop}\label{P3.10x}
 Let $F:\mathscr C\longrightarrow \mathscr D$ be a monoidal functor,  let $(U,\rho^U,\Delta,\epsilon)$ be a coalgebra object and $(R, \rho^R,\mu,\eta)$ be an algebra object in the centralizer $\mathcal Z(F)$ of $F$. Let $~\overline{\!\Theta\!}\,: U \otimes R \longrightarrow R \otimes U$ be an entwining structure in $\mathcal{Z}(F)$. Then, the bracket
\begin{equation}\label{lie}
[f,g]:=f \diamond g - (-1)^{(m-1)(n-1)} g \diamond f\qquad \qquad f\in D^m(F, U,R,~\overline{\!\Theta\!}\,),~ g\in D^n(F, U,R,~\overline{\!\Theta\!}\,) 
\end{equation} 
defines a  graded Lie algebra structure on the complex $\big(D^\bullet(F, U, R, ~\overline{\!\Theta\!}\,), \partial^\bullet\big)$. Further, the bracket $[\_\_,\_\_]$ as in \eqref{lie} and  the cup product $\cup$ as in \eqref{cup2} descend to the cohomology $H^\bullet(F, U, R, ~\overline{\!\Theta\!}\,)$ making it into a Gerstenhaber algebra.
\end{prop}
\begin{proof}
By Theorem \ref{rysz}, we know that $D^\bullet(F,U,R,~\overline{\!\Theta\!}\,)$ is a comp algebra. From the general properties of a comp algebra (see \cite{GS92}), it follows that the graded commutator in \eqref{lie} defines a  graded Lie algebra structure on the complex $D^\bullet(F, U, R,~\overline{\!\Theta\!}\,)$.  From the expression in  Lemma \ref{cupsq}, 
it is clear  that the cup product $\cup$ on $C^\bullet_{DY}(F,U,R)$  restricts to the subcomplex $D^\bullet(F,U,R,~\overline{\!\Theta\!}\,)$. Again from the general properties of a 
comp algebra in \cite{GS92}, it follows that the   cohomology $H^\bullet(F, U, R, ~\overline{\!\Theta\!}\,)$  is equipped with the bracket $[\_\_,\_\_]$ induced by \eqref{lie} and  the cup product $\cup$ induced by \eqref{cup2}, making it a Gerstenhaber algebra. 
\end{proof}

The following is now a direct consequence of the fact that  $H^\bullet(F, U, R, ~\overline{\!\Theta\!}\,)$ is a Gerstenhaber algebra. 

\begin{cor}\label{Cor3.12re}
The cup product $\cup$ on   $H^\bullet(F, U, R, ~\overline{\!\Theta\!}\,)$ is graded commutative, i.e.,
\begin{equation*}
\bar{f} \cup \bar{g}= (-1)^{mn} \bar{g} \cup \bar{f}
\end{equation*}
 for cohomology classes $\bar f \in H^m(F, U, R, ~\overline{\!\Theta\!}\,)$ and $\bar g \in \tilde H^n(F, U, R, ~\overline{\!\Theta\!}\,)$. 
\end{cor}

\begin{eg}
\emph{Suppose $\mathscr C=\mathscr D$ and let $F:\mathscr C\longrightarrow \mathscr D$ be the identity functor. Then, if $(U,\rho^U,\Delta,\epsilon)$ is a coalgebra object, $(R,\rho^R,\mu,\eta)$ is an algebra object and $~\overline{\!\Theta\!}\,: U \otimes R \longrightarrow R \otimes U$ is an entwining structure in the Drinfeld center
$\mathcal Z(id_{\mathscr C})$ of $\mathscr C$, an element $f\in D^m(F, U, R, ~\overline{\!\Theta\!}\,)$ is given by a family of morphisms   $f=\{ f_{X_1,...,X_m}: U\otimes X_1\otimes ... \otimes X_m \longrightarrow X_1\otimes ...\otimes X_m \otimes R \}_{(X_1,...,X_m)\in Ob(\mathscr C)^m}$ which satisfies the following conditions:}
\begin{align*}
(f_{X_1,...,X_m} \otimes U) \circ (U \otimes \rho^U(X_1, \ldots, X_m)) \circ (\Delta \otimes X_1 \otimes \ldots \otimes X_m)=(X_1 \otimes \ldots \otimes X_m \otimes ~\overline{\!\Theta\!}\,) \circ (\rho^U(X_1, \ldots, X_m) \otimes R) \circ 
(U \otimes f_{X_1,...,X_m}) \circ \\
(\Delta \otimes X_1 \otimes \ldots \otimes X_m)~~~~~~~~~~~~~~~~~~~~~~~~~~~~~~~~~~~~~~~~~~~~~~~~~~~~~~~~~~~~~~~~~\\
(X_1 \otimes \ldots \otimes X_m \otimes \mu) \circ (f_{X_1,...,X_m} \otimes R) \circ (U \otimes \rho^R(X_1, \ldots, X_m))=(X_1 \otimes \ldots \otimes X_m \otimes \mu) \circ (\rho^R(X_1, \ldots, X_m) \otimes R) \circ (R \otimes f_{X_1,...,X_m}) \circ \\
(~\overline{\!\Theta\!}\, \otimes X_1 \otimes \ldots \otimes X_m)~~~~~~~~~~~~~~~~~~~~~~~~~~~~~~~~~~~~~~~~~~~~~~~~~~~~~~~~~~~~~~~~~
\end{align*}
\emph{Accordingly, $D^\bullet(F, U, R, ~\overline{\!\Theta\!}\,)$ becomes a comp algebra. Further, its  cohomology groups $H^\bullet(F, U, R, ~\overline{\!\Theta\!}\,)$ are equipped with a graded commutative cup product $\cup$ and a graded Lie bracket $[\_\_,\_
\_]$ as in \eqref{lie}, giving $H^\bullet(F, U, R, ~\overline{\!\Theta\!}\,)$ the structure of a Gerstenhaber algebra.} 
\end{eg}

\begin{eg}\label{eg4.9}
\emph{Let $H$ be a finite dimensional Hopf algebra over $k$ with antipode $S$. Let $\mathscr C=\mathscr D=H-mod$, the category of left $H$-modules that are also finite dimensional over $k$. Then $H-mod$ is a finite tensor category. As explained in Example \ref{eg3.8}, it now follows that the centralizer $\mathcal Z(id_{H-mod})$ is isomorphic to the Eilenberg-Moore category $Z-mod$ of modules over the central monad $Z$ of the identity functor. From \cite[$\S$ 4.1]{Gai}, we know that this central monad $Z$ is defined as $Z(V):=H^*\otimes V$ for $V\in H-mod$, with left $H$-action given  by }
\begin{equation*}
(h\cdot (f\otimes v))(x):=f(S(h_{(1)}xh_{(3)}))\otimes h_{(2)}v \qquad f\in H^*,\textrm{ }h,x\in H, v \in V
\end{equation*}
\emph{Accordingly,  if $(U,\rho^U,\Delta,\epsilon)$ is a coalgebra object, $(R,\rho^R,\mu,\eta)$ is an algebra object and $~\overline{\!\Theta\!}\,: U \otimes R \longrightarrow R \otimes U$ is an entwining structure in $Z-mod$,  the complex $D^\bullet(F, U, R, ~\overline{\!\Theta\!}\,)$ becomes a comp algebra. Further, its  cohomology $H^\bullet(F, U, R, ~\overline{\!\Theta\!}\,)$ is equipped with a graded commutative cup product $\cup$ and a graded Lie bracket $[\_\_,\_
\_]$ as in \eqref{lie}, giving $H^\bullet(F, U, R, ~\overline{\!\Theta\!}\,)$ the structure of a Gerstenhaber algebra. A similar fact holds for the forgetful functor
$H-mod\longrightarrow Vect_k$.}

\end{eg}

\begin{eg}\label{nota}
\emph{We again consider the situation in Example \ref{eg4.9}. Let $H$ be a finite dimensional Hopf algebra. Then, it is well known that the quantum double $D(H)=H^{*op} \otimes H$ is a finite dimensional quasi-triangular Hopf algebra containing $H$ and $H^*$ as subalgebras under the association $h \mapsto 1_{H^*} \otimes h$ and $h^* \mapsto h^* \otimes 1_H$ respectively, with the universal $R$-matrix (see, for instance \cite[Theorem 7.1.1]{Majidbook} given by }
\begin{equation*}
\mathcal R=\sum \mathcal R^1 \otimes \mathcal R^2  \in H^* \otimes H \subseteq D(H) \otimes D(H)
\end{equation*}
%
\emph{Let $U$ be a coalgebra object and $R$ be an algebra object in $D(H)-mod$. 
Since $D(H)-mod \cong \mathcal Z({id_{H-mod}})$ as monoidal categories, we have  half-braidings $\rho^U$ and $\rho^R$ corresponding to $U$ and $R$ respectively,  which are given as follows: for each $X \in H-mod$, we have morphisms $\rho^U(X): U \otimes X \longrightarrow X \otimes U$ and $\rho^R(X): R \otimes X \longrightarrow X \otimes R$ in $H-mod$:
\begin{equation}\label{4.15drinfc}
\rho^U(X)(u \otimes x):= \sum \mathcal R^2 \cdot x \otimes  \mathcal R^1 \cdot u \qquad \qquad \rho^R(X)(r \otimes x):= \sum \mathcal R^2 \cdot x \otimes  \mathcal R^1 \cdot r  \qquad u \in U, r \in R, x \in X
\end{equation}
We now consider the map $~\overline{\!\Theta\!}\,: U \otimes R \longrightarrow R \otimes U$ given by
\begin{equation}\label{4.15gh}
~\overline{\!\Theta\!}\,(u \otimes r):=\sum \mathcal R^2 r \otimes  \mathcal R^1 u\qquad u\in U, r\in R
\end{equation}
where the terms $\mathcal R^2r$ and $\mathcal R^1u$ appearing in \eqref{4.15gh} are obtained from the $D(H)$-module action on $R$ and  $U$ respectively. Using the properties of the $R$-matrix $\mathcal R$, it may be verified that $~\overline{\!\Theta\!}\,$ is $D(H)$-linear, and also that $~\overline{\!\Theta\!}\,$ satisfies the conditions as in \eqref{mult} and \eqref{comult}. Therefore, $~\overline{\!\Theta\!}\,$ is an entwining structure in $D(H)-mod$. For $f \in C^m_{DY}(id_{H-mod}, U,R)$ and $X_1$,...,$X_m\in H-mod$, we write $f_{X_1,\ldots,X_m}(u \otimes x_1 \otimes \ldots \otimes x_m)=f_u^{x_1} \otimes \ldots \otimes f_u^{x_m} \otimes f_u^R \in X_1 \otimes \ldots \otimes X_m \otimes R$ for $x_i \in X_i$, $1 \leq i \leq m$ and $u \in U$. Then, a cochain in the space $D^m(id_{H-mod},U,R,~\overline{\!\Theta\!}\,)$ is an element $f \in C^m_{DY}(id_{H-mod},U,R)$ which satisfies the following two conditions:
\begin{align}
&f^{\mathcal{R}^2 \cdot {x_{1}}}_{u_{(1)}}  \otimes \ldots \otimes f^{\mathcal{R}^2 \cdot {x_{m}}}_{u_{(1)}} \otimes f^R_{u_{(1)}}  \otimes \underbrace{(\mathcal{R}^1 \ldots \mathcal{R}^1)} _{m~ \text{times}} \cdot {u_{(2)}}=\mathcal R^2 \cdot f^{x_1}_{u_{(2)}} \otimes \ldots \otimes \mathcal R^2 \cdot  f^{x_m}_{u_{(2)}} \otimes \mathcal R^2 f^R_{u_{(2)}} \otimes\underbrace{(\mathcal{R}^1  \ldots \mathcal{R}^1)} _{m+1~ \text{times}}  {u_{(1)}}\\
&f^{\mathcal R^2 \cdot x_1}_{u} \otimes \ldots \otimes f^{\mathcal R^2 \cdot x_m}_{u} \otimes \mu\big(f^R_u \otimes (\underbrace{\mathcal R^1 \ldots \mathcal R^1}_{m ~\text{times}}) \cdot r\big)=\big(\mathcal R^2 \cdot f^{x_1}_{\mathcal R^1 u}\big) \otimes \ldots \otimes \big(\mathcal R^2 \cdot f^{x_m}_{\mathcal R^1  u}\big) \otimes \mu\big(\underbrace{(\mathcal R^1  \ldots \mathcal R^1)} _{m~ \text{times}} \cdot (\mathcal R^2  r) \otimes f^{R}_{\mathcal R^1 u}\big)
\end{align}
for   $u \in U$, $r\in R$, $x_i \in X_i$, $X_i \in H-mod$,  $1 \leq i \leq m$.}
\end{eg}

\begin{eg}
\emph{Let $H$ be a finite dimensional Hopf algebra, and let $D(H)$ be the quantum double of $H$. Since $D(H)$ is also a finite dimensional Hopf algebra, its antipode $S^\circ$ is invertible (see \cite{LS}). We again take $\mathscr{C}:=H-mod$ and $F:=id_\mathscr{C}$. We now consider $D(H)$  as a left $D(H)$-module via the adjoint action, i.e.
$a \cdot_{\text{ad}} b:=a_{(1)}bS^\circ (a_{(2)})$ for all $a, b \in D(H)$, where we use   Sweedler notation for the coproduct in $D(H)$. With this adjoint action, $R:=D(H)$ becomes an algebra object in $D(H)-mod \cong \mathcal Z({id_{H-mod}})$ with the usual multiplication and unit on $D(H)$ (see \cite[$\S$ 4.1.9]{Mo}).  Similarly, it may be verified that $U:=(D(H))^*$ equipped with usual dual comultiplication and dual counit, becomes a coalgebra object in $D(H)-mod \cong \mathcal Z({id_{H-mod}})$  under the coadjoint action i.e., $(a \cdot_{\text{coad}} p)(b):=p\big(S^\circ(a_{(1)})ba_{(2)}\big)$ for all $a, b \in D(H)$ and $p \in (D(H))^*$.  Accordingly, we have half braidings $\rho^U$ and $\rho^R$ corresponding to $U=(D(H))^*$ and $R=D(H)$ respectively. For any $X \in H-mod$, $x\in X$, $p \in U$ and $a \in R$, we write $\rho^U(X)(p \otimes x)=x_\alpha \otimes p^\alpha \in X \otimes U$ and $\rho^R(X)(a \otimes x)=x_\beta \otimes a^\beta \in X \otimes R$ with summation sign omitted. 
}

\smallskip
\emph{We now consider the map $~\overline{\!\Theta\!}\,: (D(H))^* \otimes D(H) \longrightarrow D(H) \otimes (D(H))^*$ defined by
\begin{equation*}
~\overline{\!\Theta\!}\,(p \otimes a):=a_{(2)} \otimes {S^\circ}^{-1}(a_{(1)}) \cdot_{\text{coad}} p \qquad p \in D(H))^*,~ a \in D(H) 
\end{equation*}
It may be verified that the map $~\overline{\!\Theta\!}\,$ is an entwining structure in the category $D(H)-mod$.  Using the same notation as in Example \ref{nota}, a cochain in the   subspace $D^m(id_{H-mod}, U, R,~\overline{\!\Theta\!}\,)$ is a morphism $f \in C^m_{DY}(id_{H-mod}, U, R)$ which satisfies the following two conditions:
\begin{align}
 f^{{x_1}_{\alpha_1}}_{p_{(1)}} \otimes \ldots  \otimes f^{{x_m}_{\alpha_m}}_{p_{(1)}} \otimes f^R_{p_{(1)}} \otimes (p_{(2)})^{\alpha_1 \cdots \alpha_m}=\big(f^{x_1}_{p_{(2)}}\big)_{\alpha_1} \otimes \ldots \otimes \big(f^{x_m}_{p_{(2)}}\big)_{\alpha_m} \otimes \big(f^{R}_{p_{(2)}}\big)_{(2)} \otimes {S^\circ}^{-1}\big(\big(f^{R}_{p_{(2)}}\big)_{(1)}\big) \cdot_{\text{coad}} (p_{(1)})^{\alpha_1\cdots \alpha_m}\\
f_p^{{x_1}_{\beta_1}} \otimes \ldots \otimes f_p^{{x_m}_{\beta_m}} \otimes \mu(f_p^R \otimes  a^{\beta_1 \cdots \beta_m}) = \big(f^{x_1}_{{S^\circ}^{-1}(a_{(1)}) ~\cdot_{\text{coad}}~ p} \big)_{\beta_1} \otimes \ldots \otimes \big(f^{x_m}_{{S^\circ}^{-1}(a_{(1)}) ~\cdot_{\text{coad}}~ p} \big)_{\beta_m} \otimes  \mu\big(a_{(2)})^{\beta_1\cdots \beta_m} \otimes  f^{R}_{{S^\circ}^{-1}(a_{(1)}) ~\cdot_{\text{coad}}~ p}\big)
\end{align}
for $X_1$, ..., $X_m\in H-mod$ and $x_i\in X_i$ for $1\leq i\leq m$.}
\end{eg}

Our final aim in this paper is to dualize the constructions in the first part of this section. This will lead to another subcomplex of $C^\bullet_{DY}(F,U,R)$ corresponding to a comp algebra, and whose cohomology is therefore a Gerstenhaber algebra. For this, we consider a morphism $~\underline{\!\Psi\!}\,: R \otimes U \longrightarrow U \otimes R$ in the category $\mathcal{Z}(F)$ which is also compatible with the algebra structure of $R$. This means that the following diagrams commute:

\begin{equation}\label{morpxx}
\begin{tikzcd}[column sep=large, row sep=large]
R \otimes U \otimes F(X)
  \arrow[r, "R \otimes \rho^U(X)"]
  \arrow[d, swap,"~\underline{\!\Psi\!}\, \otimes {F(X)}"]
&
R \otimes F(X) \otimes U
  \arrow[r, "\rho^R(X) \otimes U"]
&
F(X) \otimes R \otimes U
  \arrow[d, "{F(X)} \otimes ~\underline{\!\Psi\!}\,"]
\\
U \otimes R \otimes F(X)
  \arrow[r, "U \otimes \rho^R(X)"]
&
U \otimes F(X) \otimes R
  \arrow[r, "\rho^U(X) \otimes R"]
&
F(X) \otimes U \otimes R
\end{tikzcd}
\end{equation}
\begin{equation}\label{oppent1}
\begin{tikzcd}
 R \otimes R \otimes U
\arrow[r, "{R \otimes ~\underline{\!\Psi\!}\,}"]
\arrow[d, swap, "{\mu \otimes U}"]
&
R \otimes U \otimes R
\arrow[r, "{~\underline{\!\Psi\!}\, \otimes R}"]
&
U \otimes R \otimes R
\arrow[d, "{U \otimes \mu}"]\\
R \otimes U
\arrow[rr,"{~\underline{\!\Psi\!}\,}"]
& 
&
U \otimes R 
\end{tikzcd} \qquad \qquad
\begin{tikzcd}
& U 
\arrow[ld, swap, "{\eta \otimes U }"]
\arrow[rd, "{U \otimes \eta}"]
&\\
R \otimes U
\arrow[rr, "{~\underline{\!\Psi\!}\,}"] 
& & U \otimes R
\end{tikzcd}
\end{equation}

For each $m \geq 0$, we define $P^m_{DY}(F,U,R, ~\underline{\!\Psi\!}\,)\subseteq C^m_{DY}(F,U,R)=Nat(U\otimes F^{\otimes m},F^{\otimes m}\otimes R)$ to be the collection of those natural tranformations $f:U \otimes F^{\otimes m} \longrightarrow F^{\otimes m} \otimes R$ for which the following diagram commutes
\begin{equation}\small
\begin{tikzcd}[column sep=2.5cm, row sep=large]
U \otimes F(X_1) \otimes \ldots \otimes F(X_m)
\arrow[r, "{\Delta \otimes F(X_1) \otimes \ldots \otimes F(X_m)}"]
&
U \otimes U \otimes F(X_1) \otimes \ldots \otimes F(X_m)
\arrow[r, "{U \otimes f_{X_1,\ldots,X_m}}"]
\arrow[d, swap, "{U \otimes \rho^U(X_1,\ldots,X_m)}"]
&
U \otimes F(X_1) \otimes \ldots \otimes F(X_m) \otimes R
\arrow[dd,"{\rho^U(X_1,\ldots,X_m) \otimes R}"]
\\
& U \otimes F(X_1) \otimes \ldots \otimes F(X_m) \otimes U
\arrow[d, swap, "{f_{X_1,\ldots,X_m}  \otimes U}"]
\\
&
F(X_1) \otimes \ldots \otimes F(X_m) \otimes R \otimes U
\arrow[r, "{F(X_1) \otimes \ldots \otimes F(X_m) \otimes ~\underline{\!\Psi\!}\,}"]
&
F(X_1) \otimes \ldots \otimes F(X_m) \otimes U \otimes R
 \end{tikzcd}
 \end{equation}
for all $X_1$, ..., $X_m\in Ob(\mathscr C)$. Similar to Proposition \ref{ent11} and Lemma \ref{closed}, this gives us the following result.

\begin{prop}
The subspaces $P^m(F,U,R, ~\underline{\!\Psi\!}\,)\subseteq C^m_{DY}(F,U,R)$ for $m\geq 0$ determine a subcomplex of $(C^\bullet_{DY}(F,U,R),\partial^\bullet)$. Further,  $f \diamond_i g \in P^{m+n-1}(F,U,R,~\underline{\!\Psi\!}\,)$ for any $f \in P^{m}(F,U,R,~\underline{\!\Psi\!}\,)$ and $g \in P^{n}(F,U,R,~\underline{\!\Psi\!}\,)$ for $m$, $n\geq 0$.
\end{prop}

\smallskip
Next, we consider a morphism $~\underline{\!\Phi\!}\,: R \otimes U \longrightarrow U \otimes R$ in the category $\mathcal{Z}(F)$ compatible with the coalgebra structure of $U$. This means that, in addition to the condition as in \eqref{morpxx}, the following diagrams commute:
\begin{equation}\label{oppent2}
\begin{tikzcd}
R \otimes U
\arrow[r, "{R \otimes \Delta}"]
\arrow[d, swap, "{~\underline{\!\Phi\!}\,}"]
&
R \otimes U \otimes U 
\arrow[r, "{~\underline{\!\Phi\!}\, \otimes U}"]
&
U \otimes R \otimes U
\arrow[d, "{U \otimes ~\underline{\!\Phi\!}\,}"]\\
U \otimes R
\arrow[rr,"{\Delta \otimes R }"]
& 
&
 U \otimes U \otimes R
\end{tikzcd}
\qquad \qquad
\begin{tikzcd}
R \otimes U
\arrow[rr,"{~\underline{\!\Phi\!}\,}"]
\arrow[dr, swap, "{R \otimes \epsilon }"]
&
&
U \otimes R
\arrow[dl, "{\epsilon \otimes R}"]\\
& R
\end{tikzcd}
\end{equation}
Now, for each $m \geq 0$, we define $Q^m(F,U,R, ~\underline{\!\Phi\!}\,)\subseteq C^m_{DY}(F,U,R)=Nat(U\otimes F^{\otimes m},F^{\otimes m}\otimes R)$ to be the collection of those natural tranformations $f:U \otimes F^{\otimes m} \longrightarrow F^{\otimes m} \otimes R$ for which the following diagram commutes
\begin{equation}\small
\begin{tikzcd}[column sep=2.5cm, row sep=large]
R \otimes U \otimes F(X_1) \otimes \ldots \otimes F(X_m)
\arrow[r, "{R \otimes f_{X_1,\ldots,X_m}}"]
\arrow[d, "{~\underline{\!\Phi\!}\, \otimes F(X_1) \otimes \ldots \otimes F(X_m)}"']
&
R \otimes F(X_1) \otimes \ldots \otimes F(X_m) \otimes R
\arrow[dd, "{\rho^R(X_1,\ldots,X_m) \otimes R}"]
&
\\
U \otimes R \otimes F(X_1) \otimes \ldots \otimes F(X_m)
\arrow[d, swap, "{U \otimes \rho^R(X_1,\ldots,X_m)}"]
\\
U \otimes F(X_1) \otimes \ldots \otimes F(X_m) \otimes R
\arrow[r, "{f_{X_1,\ldots,X_m} \otimes R}"]
&
F(X_1) \otimes \ldots \otimes F(X_m) \otimes R \otimes R
\arrow[r, "{F(X_1) \otimes \ldots \otimes F(X_m) \otimes \mu}"]
&
F(X_1) \otimes \ldots \otimes F(X_m) \otimes R 
 \end{tikzcd}
 \end{equation}
for all $X_1$, ..., $X_m\in Ob(\mathscr C)$. Similar to Proposition \ref{prop4.3} and Lemma \ref{closedn}, this gives us the following result.

\begin{prop}
The subspaces $Q^m(F,U,R,~\underline{\!\Phi\!}\,)\subseteq C^m_{DY}(F,U,R)$ for $m\geq 0$ determine a subcomplex of $(C^\bullet_{DY}(F,U,R),\partial^\bullet)$. Further,  $f \diamond_i g \in Q^{m+n-1}(F,U,R,~\underline{\!\Phi\!}\,)$ for any $f \in Q^{m}(F,U,R,~\underline{\!\Phi\!}\,)$ and $g \in Q^{n}(F,U,R,~\underline{\!\Phi\!}\,)$ for $m$, $n\geq 0$.
\end{prop}

\smallskip
We now let $~\underline{\!\Theta\!}\,: R \otimes U \longrightarrow U \otimes R$ be an entwining structure in the category $\mathcal{Z}(F)$, i.e., $~\underline{\!\Theta\!}\,$ is a morphism in $\mathcal{Z}(F)$ satisfying both the conditions as in \eqref{oppent1} and \eqref{oppent2}. For each $m \geq 0$, we set
\begin{equation}
D^m(F,U,R,~\underline{\!\Theta\!}\,):= P^m(F,U,R,~\underline{\!\Theta\!}\,) \cap Q^m(F,U,R,~\underline{\!\Theta\!}\,)\subseteq C^m_{DY}(F,U,R)
\end{equation}
It is   clear that $\big(D^\bullet(F,U,R,~\underline{\!\Theta\!}\,), \partial^\bullet\big)$ is a subcomplex of $\big(C^\bullet_{DY}(F,U,R), \partial^\bullet\big)$. We conclude with the following result.

\begin{Thm}\label{Th4.14ed}
 Let $F:\mathscr C\longrightarrow \mathscr D$ be a monoidal  functor between monoidal  categories, let $(U,\rho^U,\Delta,\epsilon)$ be a coalgebra object and $(R, \rho^R,\mu,\eta)$ be an algebra object in the centralizer $\mathcal Z(F)$ of $F$. Let $~\underline{\!\Theta\!}\,: R \otimes U \longrightarrow U \otimes R$ be an entwining structure in $\mathcal{Z}(F)$. Then, 

\smallskip
(a)  $P^\bullet(F,U,R,~\underline{\!\Theta\!}\,) $ and $Q^\bullet(F,U,R,~\underline{\!\Theta\!}\,) $  are weak comp algebras.

\smallskip 
(b)  $ D^\bullet(F,U,R,~\underline{\!\Theta\!}\,) $ is a comp algebra.
 
\smallskip
(c)  The bracket
\begin{equation}\label{liembsk}
 [f,g]:=f \diamond g - (-1)^{(m-1)(n-1)} g \diamond f\qquad \qquad f\in D^m(F, U,R,~\underline{\!\Theta\!}\,),~ g\in D^n(F, U,R,~\underline{\!\Theta\!}\,) 
 \end{equation}
defines a  graded Lie algebra structure on the complex  $\big(D^\bullet(F,U,R,~\underline{\!\Theta\!}\,), \partial^\bullet\big)$. Further, the bracket $[\_\_,\_\_]$ as in \eqref{liembsk} and  the cup product $\cup$ as in \eqref{cup2} descend to  the cohomology $H^\bullet(F,U,R,~\underline{\!\Theta\!}\,)$ of the  complex $\big(D^\bullet(F,U,R,~\underline{\!\Theta\!}\,), \partial^\bullet\big)$, making $H^\bullet(F,U,R,~\underline{\!\Theta\!}\,)$  into a Gerstenhaber algebra.   In particular, the cup product $\cup$ on $H^\bullet(F,U,R,~\underline{\!\Theta\!}\,)$ is graded commutative.
\end{Thm}
\begin{proof}
This may be proved in a manner similar to Theorem \ref{rysz}, Proposition \ref{P3.10x} and Corollary \ref{Cor3.12re}.
\end{proof}

\end{document}